\documentclass[12pt]{elsarticle}
\usepackage[margin=1.5cm]{geometry}
\usepackage{amsmath}
\usepackage{amssymb}
\usepackage{mathrsfs}
\usepackage{dsfont}
\usepackage{amsthm}
\usepackage{cases}
\usepackage{hyperref}
\usepackage{abstract}
\usepackage{float}
\usepackage{extarrows}
\usepackage[dvipsnames]{xcolor}
\usepackage{enumerate}
\usepackage[shortlabels]{enumitem}

\journal{a journal}
\date{}

\hypersetup{
    colorlinks=True,
    linkcolor=blue,
    filecolor=magenta,      
    urlcolor=cyan,
    citecolor=ForestGreen
}
\newtheorem{Theorem}{Theorem}[section]
\newtheorem{Proposition}{Proposition}[section]

\newtheorem{Corollary}{Corollary}[section]

\theoremstyle{definition}

\newtheorem{Definition}{Definition}[section]
\newtheorem{Remark}{Remark}[section]

\def\ra{\rangle}
\def\la{\langle}
\def\intOm{\int\limits_\Omega}

\begin{document}



\title{Stress solution of static linear elasticity with mixed boundary conditions via adjoint linear operators} 


\author[1]{Ivan Gudoshnikov\corref{cor1} \fnref{fn2} }
\author[1]{Michal K\v{r}\'{i}\v{z}ek \fnref{fn3}}
\cortext[cor1]{Corresponding author}
\fntext[fn2]{gudoshnikov@math.cas.cz, ORCiD: 0000-0002-6136-484X}
\fntext[fn3]{krizek@math.cas.cz, ORCiD: 0000-0002-5271-5038}
\affiliation[1]{organization={Institute of Mathematics of the Czech Academy of Sciences},
            addressline={\v{Z}itn\'{a} 25}, 
            city={Praha 1},
            postcode={115 67}, 
            country={Czech Republic}}
\maketitle
\begin{abstract}
We revisit stress problems in linear elasticity to provide a perspective from the geometrical and functional-analytic points of view. For the static stress problem of  linear elasticity with mixed boundary conditions we present the associated pair of unbounded adjoint operators. Such a pair is explicitly written for the first time, despite the abundance of the literature on the topic. We use it to find the stress solution as an intersection of the (affinely translated) fundamental subspaces of the adjoint operators. In particular, we treat  the equilibrium equation in the operator form, which involves the spaces of traces on a part of the boundary, known as the Lions-Magenes spaces. Our analysis of the pair of adjoint operators for the problem with mixed boundary conditions relies on the properties of the analogous pair of operators for the problem with the displacement boundary conditions, which we also include in the paper.
\end{abstract}






\section{Introduction} 
Static linear elasticity is the cornerstone of continuum mechanics as a scientific, engineering and mathematical discipline, and its development spans from the beginning of the nineteenth century to the current day, see \cite{Truesdell1984} for historical remarks.  The mathematical problems involved are elliptic partial differential equations with vector-valued or matrix-valued distributions as the unknown variables. To find a weak solution, such problems are usually rewritten in a variational form, which suggests a minimization of a convex lower semicontinuous functional (see e.g. \cite[Section 7.2]{Necas1981}, \cite[Section 4.6]{Sayas2019}). Such a form is also convenient for defining finite element  approximations (see many books, including \cite{Boffi2013}, \cite{Ciarlet1978}, \cite{Ern2004}, \cite{Gatica2014}, \cite{Krizek1990}, which pay a particular attention to elasticity). Additionally, for the elasticity problem it is possible to write a self-adjoint second-order differential operator, analogous to the scalar Laplace operator \cite[Section 10]{McLean2000}, which enables one to apply spectral theory.

However, to see the geometric structure of the problem it is useful to employ the method of orthogonal subspaces, similar to Beltrami decomposition (see \cite[(5.15)]{Admal2016}, cf. \cite{Geymonat2005}), but taking into account boundary conditions (see e.g. \cite[Section 5.1, Th. 5.1]{Amstutz2020}). 

In the paper we focus on the method of orthogonal subspaces applied to the static problem of  linear elasticity with mixed boundary conditions. We focus on the stress variable as the main unknown, cf. \cite{Ciarlet2004}, where the main unknown is the strain variable and the boundary condition is pure traction. Our goal is to rigorously describe the stress solution in terms of the orthogonal subspaces. Such description will be instrumental for the analysis of the models with more complex constitutive relations involving stress variables, e.g. with plasticity \cite{Han2012} or piezoelectricity \cite{Kaltenbacher2007}.

\subsection*{The method of orthogonal subspaces}
\begin{figure}[H]\center
\includegraphics{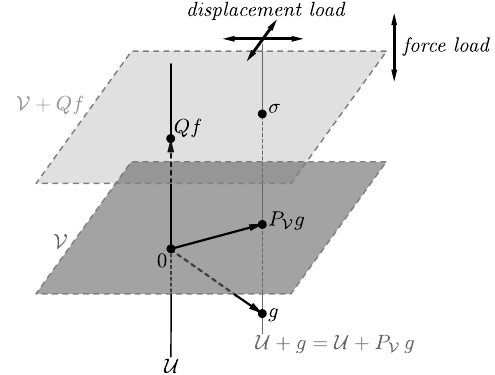}
\caption{\footnotesize Schematic representation of static stress problems of linear elasticity. The subspace $\mathcal{U}$ (depicted vertically) represents the stresses, which correspond, in the sense of Hooke's law, to the image of the ``strain'' operator (${\bf C}\,{\rm Im}\, {\bf E_0}$ or ${\bf C}\,{\rm Im}\, {\bf E_1}$ in terms of the paper). The subspace $\mathcal{V}$ (depicted horizontally) represents the kernel of the ``divergence'' operator (respectively, ${\rm Ker}\, {\bf D_0}$ or ${\rm Ker}\, {\bf D_1}$). The subspaces $\mathcal{U}$ and $\mathcal{V}$ are orthogonal complements in the Hilbert space of stress distributions, and the operator $P_{\,\mathcal{V}}$ is the orthogonal projection onto $\mathcal{V}$. Affine sets $\mathcal{U}+g$ and $\mathcal{V}+Qf$ are shifted according to prescribed boundary displacements and prescribed force loads, respectively. A stress solution $\sigma$ is located at the intersection of the orthogonal affine sets, therefore it is unique.
} 
\label{fig:spaces}
\end{figure}

The main idea of the method of orthogonal subspaces comes from the fact that the unknown stress variable must satisfy two groups of equations. The first group is the kinematic compatibility of the corresponding strain and the displacement boundary conditions, while the second group consists of the equilibrium equation. These equations can be written in terms of the unknown stress  variable from the space $L^2_{\rm sym}(\Omega)$ of square-integrable fields with values in symmetric matrices (``symmetric matrix-valued fields'') and defined over a bounded connected domain $\Omega$ with Lipschitz boundary. The two groups of equations mean that the unknown stress variable must lay at the intersection of two orthogonal affine sets (see Fig. \ref{fig:spaces}), which are, up to parallel translations, 
\begin{itemize}
\item 
$\mathcal{U}$ --- the subspace of stresses, whose corresponding strains satisfy kinematic compatibility and homogeneous displacement boundary conditions,
\item 
$\mathcal{V}$ --- the subspace of divergence-free stresses with zero trace on the free boundary (if it is present), i.e. the stress distributions with vanishing resultant forces.
\end{itemize}
These subspaces are orthogonal complements in $L^2_{\rm sym}(\Omega)$ (with the inner product weighted according to the elasticity tensor, to be precise). In turn, the parallel translations  depend linearly on prescribed displacements on the boundary and prescribed force loads.

The observation above provides a simple and intuitive way to think about static linear elasticity. Indeed, despite many complex details, finding the stress distribution is nothing but finding an intersection of two closed affine sets orthogonal to each other, see Fig. \ref{fig:spaces} again. Notice, for example, that the uniqueness follows trivially from such a description.

\subsection*{The adjoint operators}

The important fact, which is rarely discussed however, is that the orthogonality of the subspaces follows from the adjacency between the operators of strain and divergence, and  their modifications, which take into account various boundary conditions. In the case of continuum mechanics the operators are defined with the values in symmetric matrix-valued fields and vector fields, respectively. The adjacency between such operators is occasionally used in the literature (e.g. \cite[the proof of Prop.~4.4, pp.~3779--3780]{Amstutz2020}, \cite[the proof of Th. 5, p.~177]{Geymonat2005}, \cite[p.~210]{Geymonat1986}), but the only dedicated study we found was about the simpler case, when the gradient and the divergence operators take values in vector fields and scalar fields, respectively \cite{Kurula2012}.

Our main objective is to solve the stress problem in static linear elasticity with mixed boundary conditions (displacement prescribed on the part $\Gamma_{\rm D}$, traction force prescribed on the part $\Gamma_{\rm T}$) using the tools of linear functional analysis, and the theory of unbounded adjoint operators in particular. This will allow us to understand the geometric structure of the problem (Fig. \ref{fig:spaces}), which can be overlooked when a general variational framework is used. 

The second objective is to explicitly consider the mapping between stress distributions and resultant forces in the body volume and on the part of the boundary $\Gamma_{\rm T}$). While it is possible to solve the linear elasticity problem without an explicit consideration of the stress-to-forces mapping (see e.g. \cite[Chapter 7]{Necas1981}), the natural questions are:
\begin{itemize}
\item What is the proper way to set up such a mapping and what are the spaces of functions involved?
\item What are the properties of the mapping, in particular, regarding its surjectivity and closedness of the associated spaces?
\end{itemize}
The answer to the first question above is to write the stress-to-forces mapping as
\begin{gather*}
{\bf D_1}: H_{\rm sym}({\rm div}, \Omega) \subset  L^2_{\rm sym}(\Omega)\to [L^2(\Omega)]^n \times [H^{-1/2}_{00}(\Gamma_{\rm T})]^n,
\\
{\bf D_1}: \sigma\mapsto \begin{pmatrix} - {\rm div}\, \sigma \\  {\rm tr}_{\nu, \Gamma_{\rm T}} (\sigma)
\end{pmatrix},
\end{gather*}
where $H_{\rm sym}({\rm div}, \Omega)$ is the space of symmetric matrix-valued functions which admit weak divergence, ${\rm tr}_{\nu, \Gamma_{\rm T}}$ is the normal trace on the part of the boundary, and the space $H^{-1/2}_{00}(\Gamma_{\rm T})$ is the dual Lions-Magenes space (short overviews will be provided in Sections \ref{ssect:weak_divergence}, \ref{ssect:L-M_spaces}). Here and in the entire text $n\geqslant 2$ is the number of physical dimensions, so that $\Omega \subset \mathbb{R}^n$. 

In the paper we systematically follow the theory of unbounded adjoint operators to derive all the properties of ${\bf D_1}$ and its adjoint operator, which are necessary to solve the elasticity problem. The central theorem of the text (Theorem \ref{th:main-pair}) states that ${\bf D_1}$ and the operator 
\begin{gather*}
{\bf E_1}: D({\bf E_1}) = \left\{(u, {\rm tr\,}_{\Gamma_{\rm T}}(u)): u\in [W_{\rm D}]^n\right\}\subset [L^2(\Omega)]^n\times [H^{1/2}_{00}(\Gamma_{\rm T})]^n \to L^2_{\rm sym}(\Omega),
\\
{\bf E_1}:\begin{pmatrix} u\\ {\rm tr\,}_{\Gamma_{\rm T}}(u)\end{pmatrix} \mapsto \frac{1}{2}\left(\nabla u+ (\nabla u)^{\top}\right),
\end{gather*}
are mutually adjoint, and that both ${\bf D_1}, {\bf E_1}$ are closed and densely defined ($W_{\rm D}$ above is the space of $W^{1,2}$-functions, vanishing on $\Gamma_{\rm D}$). The properties of adjoint operators allow us to establish the surjectivity of ${\bf D_1}$ from Korn's second inequality. The solvability of the elasticity problem and the applicability of Fig. \ref{fig:spaces} then follow for any applied body force and traction force loads.

\subsection*{A link between the problem with mixed boundary conditions and the problem with the displacement boundary conditions.}

We would like to point out a particular obstacle in the proof of the fact that ${\bf E_1}={\bf D_1}^*$, which we discussed above. Observe, that if one considers an arbitrary element of the dual space to the range of ${\bf D_1}$,
\begin{equation}
\left(u, {\rm tr\,}_{\Gamma_{\rm T}}(v)\right)\in \left([L^2(\Omega)]^n \times [H^{-1/2}_{00}(\Gamma_{\rm T})]^n\right)^* \cong [L^2(\Omega)]^n \times [H^{1/2}_{00}(\Gamma_{\rm T})]^n,
\label{eq:intro-about-bad-functional}
\end{equation}
its components are two {\it different} functions, of which $v$ is a Sobolev function from $[W_{\rm D}]^n$ and $u$ is merely from $[L^{2}(\Omega)]^n$. Thus, in general, such an element $\left(u, {\rm tr\,}_{\Gamma_{\rm T}}(v)\right)$ has much worse regularity than the elements of $D({\bf E_1})$.

Interestingly, the way to overcome this obstacle highlights the connection between the problem with mixed boundary conditions and the problem with the displacement boundary conditions only (i.e. when $\Gamma_{\rm D}$ is the entire boundary). Specifically, from the definition of the domain of the adjoint operator and Green's formula it follows (see \eqref{eq:D(-Div)-def} below and the derivation afterwards) that 
\begin{equation}
u-v\in D({\bf D_0}^*),
\label{eq:passage-from-D1-to-D0}
\end{equation}
where ${\bf D_0}$ is the stress-to-forces operator for the problem with the displacement boundary conditions, written as follows:
\begin{gather*}
{\bf D_0}:D({\bf D_0}) = H_{\rm sym}({\rm div}, \Omega) \subset  L^2_{\rm sym}(\Omega)\to [L^2(\Omega)]^n,
\\
{\bf D_0}: v\mapsto - {\rm div}\, v.
\end{gather*}
To be able to obtain the regularity of $u$ from \eqref{eq:passage-from-D1-to-D0} we first prove Theorem \ref{th:auxiliary_operators_are_mutually_adjoint}, which states that ${\bf D_0}$ and the operator
\begin{gather*}
{\bf E_0}: D({\bf E_0}) = [W_0^{1,2}(\Omega)]^n \subset [L^2(\Omega)]^n \to L^2_{\rm sym}(\Omega), \\
{\bf E_0}:u \mapsto \frac{1}{2}\left(\nabla u+ (\nabla u)^{\top}\right),
\end{gather*}
are mutually adjoint, closed and densely defined, in the similar manner to ${\bf D_1}$ and ${\bf E_1}$. Then \eqref{eq:passage-from-D1-to-D0} yields that
\[
u-v \in [W_0^{1,2}(\Omega)]^n,
\]
i.e. $u$ is a Sobolev function with the same trace as $v$. We conclude that, for an element $\left(u, {\rm tr\,}_{\Gamma_{\rm T}}(v)\right)$ satisfying \eqref{eq:intro-about-bad-functional}, if it is in $D({\bf D_1}^*)$, then it must belong to $D({\bf E_1})$.

This way the operators ${\bf E_0}, {\bf D_0}$ (which describe the elasticity problem with displacement boundary conditions) are used to establish the properties of the operators ${\bf E_1}, {\bf D_1}$ (which describe the elasticity problem with mixed boundary conditions). Additionally, as a byproduct of Theorem \ref{th:auxiliary_operators_are_mutually_adjoint} about ${\bf D_0}$ and ${\bf E_0}$, we solve the elasticity problem with displacement boundary conditions.

\subsection*{A connection to discrete models}
It is worth to note, that the framework of the adjoint operators and the associated orthogonal subspaces (Fig. \ref{fig:spaces}) is common for spatially continuous and spatially discrete models (\cite[Th. 1, pp.~79--80]{Admal2016}, \cite[Section 3]{Gudoshnikov2023preprint}, \cite[(3.5)--(3.7), p.~7]{Gudoshnikov2021}), because it is the consequence of the {\it principle of virtual work}, general for mechanics \cite[Section 1.4, pp.~16--17]{Goldstein2002}.

For the discrete models, the role of the unbounded adjoint operators of strain and divergence is taken by their matrix analogues, the {\it compatibility matrix} and the {\it equilibrium matrix}  \cite{Lubensky2015},  \cite[Section 3]{Gudoshnikov2023preprint}. The properties of such matrices and the properties of the underlying graph (``bar framework'') such as {\it rigidiy} \cite{RothWhiteley1981}, \cite[Chapters 8 and 9]{Alfakih2018_rigidity}, are important for the elasticity problem in the same manner as the properties of the strain operator (e.g. Proposition \ref{prop:rigid_motions} below) is important for the elasticity problem of continuum media. In the particular case of a discrete model with one spatial dimension, the discrete analogue of the divergence operator is the incidence matrix of the underlying directed graph \cite[Chapters 2 and 5]{Bapat2010}, \cite[Section 2]{Gudoshnikov2021}.

For all such discrete elasticity models the Fundamental Theorem of Linear Algebra (see e.g. \cite[Th. 4.45, p.~221]{OlverLinearAlgebra2018}) yields the picture, similar to Fig. \ref{fig:spaces}, and  allows to find a unique stress solution. With the current text we show how the intuition about the transpose of a matrix and related orthogonality of the kernels and the images in the Fundamental Theorem of Linear Algebra can be carried over to problems in continuum mechanics. We also hope that our text could potentially serve as an introduction to linear elasticity for the people with a background in functional analysis.

\subsection*{The contents of the text}

The paper is organized as follows. In Section \ref{sect:general_preliminaries} we remind the reader the fundamental mathematical preliminaries. There we also include the framework of unbounded adjoint operators, which plays the major role for us, and a related discussion on solving an abstract linear system written via an unbounded operator. In Section \ref{sect:datum} we formulate the linear elasticity problems in terms of the basic function spaces, unknown variables, given data and the equations involved. Section \ref{sect:problem1} starts with a brief description of the space $H_{\rm sym}({\rm div}, \Omega)$, which is then used to construct the pair of unbounded adjoint operators ${\bf E_0}$, ${\bf D_0}$ and to solve the stress problem with the displacement boundary conditions. Section \ref{sect:problem2} begins with an introduction of the traces on a part of the boundary and the associated function spaces (Lions-Magenes spaces). Using such spaces we construct the pair of unbounded adjoint operators ${\bf E_1}$, ${\bf D_1}$ to solve the stress problem with the mixed boundary conditions. Pair ${\bf E_0}$, ${\bf D_0}$, to which the previous section is largely devoted, appears useful here in the proof of the adjacency between ${\bf E_1}$ and ${\bf D_1}$. Section \ref{sect:conclusions} contains some concluding remarks, and in 
\ref{sect:appendix-on-other-problems} we discuss a connection to the well-known variational formulations and the displacement problems using the notation of the current paper. Apart from the Appendix, we do not pursue the solution in terms of the displacement variable, and focus on the stress problem. In particular, we do not discuss Saint-Venant's representations of the involved spaces via operator ${\rm curl\, curl}$, for which we refer to e.g. \cite{Geymonat2005}, \cite[Th. 31, p.~37]{Necas1981}.  We must note, that static linear elasticity and similar problems can also be examined in terms of de Rham Hilbert complexes, for which we refer to \cite{PauliSchomburg2022} and references therein.

\section{General preliminaries}
\subsection{Notation}
For a linear operator $A$ we use the notation ${\rm Ker}\, A$ and ${\rm Im}\,A$ for the kernel (the nullspace) and the image of $A$, respectively. When $A$ is bounded, we denote by $\|A\|_{\rm op}$ the operator norm of $A$ (see e.g. \cite[p.~43]{Brezis2011}).
\label{sect:general_preliminaries}
\subsection{Domain}
Throughout the text $n\geqslant 2$ and {\it the domain} $\Omega\subset\mathbb{R}^n$ is an open bounded connected nonempty set with Lipschitz boundary $\Gamma$. For the definition and discussion on Lipschitz boundary see e.g. \cite[Section 1.2]{Grisvard2011}. In particular, the Lipschitz property of the boundary implies the existence of a unit outward normal vector 
\[
\nu(x) = (\nu_1(x),\dots, \nu_n(x))\in \mathbb{R}^n
\]
for a.a. $x\in \Gamma$, see \cite[p.~37]{Grisvard2011}. Throughout the text the integration over $\Omega$ is meant with respect to the Lebesgue measure and the integration over $\Gamma$ and its open (in $\Gamma$) subsets is with respect to the Hausdorff measure $\mathcal{H}^{n-1}$, see e.g. \cite[Appendix C]{Leoni2017}.
\subsection{Tensor calculus}
We would like to remind the following basic definitions.
\begin{Definition} For a function $u:\Omega\subset \mathbb{R}^n \to \mathbb{R}^n$ the {\it gradient} is defined as
\[
\nabla u =\begin{pmatrix} \frac{\partial u_1}{\partial x_1} &\cdots &\frac{\partial u_1}{\partial x_n} \\ \vdots & & \vdots\\
\frac{\partial u_n}{\partial x_1} &\cdots &\frac{\partial u_n}{\partial x_n}
\end{pmatrix},
\]
i.e. when $u$ is differentiable in the classical sense, $\nabla u$ coincides with the Jacobi matrix (Fr\'{e}chet derivative). In turn, let $\sigma:\Omega\subset \mathbb{R}^n \to \mathbb{R}^{n\times n}$ be a 2-tensor field. Then its divergence in the classical sense is defined as
\[
{\rm div}\, \sigma = \nabla \cdot \sigma  = 
\begin{pmatrix}
\frac{\partial \sigma_{11}}{\partial x_1}+ \cdots + \frac{\partial \sigma_{1n}}{\partial x_n}  \\ \vdots\\
\frac{\partial \sigma_{n1}}{\partial x_1}+ \cdots + \frac{\partial \sigma_{nn}}{\partial x_n}
\end{pmatrix}.
\]
 \end{Definition}
\subsection{Sobolev spaces and traces}
The introduction to Sobolev spaces is available in a plethora of textbooks, see e.g. \cite[Chapter 3]{McLean2000}, \cite[Chapters 10--11]{Leoni2017}, \cite[Chapter 5]{Evans2010}. We only remind that $W^{1,2}(\Omega)$ is a Hilbert space, in which the space $\mathcal{D}\left(\overline{\Omega}\right)$ of infinitely differentiable functions is dense, see \cite[Th. 11.35, p.]{Leoni2017}, \cite[Th. 1.4.2.1, p.~24]{Grisvard2011}. To avoid a confusion about notation, we also refer to the standard definitions of $\mathcal{D}({\Omega})$ and $\mathcal{D}\left(\overline{\Omega}\right)$, see  \cite[pp.~61, 65,77]{McLean2000}.

We will use the following {\it intrinsic} definition of a  fractional Sobolev space.
\begin{Definition}
For $\Gamma'=\Gamma$ or $\Gamma'$ being an open in $\Gamma$ subset of dimension $n-1$ denote
\begin{equation}
\la u, v \ra_{H^{1/2}({\Gamma'})} = \int\limits_{\Gamma'} u\, v\, dS + \int\limits_{\Gamma'} \int\limits_{\Gamma'} \frac{(u(x)-u(y))(v(x)-v(y))}{\|x-y\|^n}dS_x dS_y,
\label{eq:H1/2-ip-def}
\end{equation}
\begin{equation}
\|u\|_{H^{1/2}({\Gamma'})} = \left(\la u,u\ra_{H^{1/2}({\Gamma'})} \right)^{\frac{1}{2}} = \left(\|u\|_{L^2(\Gamma')}^2+ \int\limits_{\Gamma'} \int\limits_{\Gamma'} \frac{|u(x)-u(y)|^2}{\|x-y\|^n}dS_x dS_y\right)^{\frac{1}{2}},
\label{eq:H1/2-norm-def}
\end{equation}
and define the the {\it fractional Sobolev space}
\begin{equation*}
H^{1/2}(\Gamma') = \left\{ u\in L^2(\Gamma'): \|u\|_{H^{1/2}({\Gamma'})} <\infty  \right\}.
\end{equation*}
\end{Definition}
Note that $H^{1/2}({\Gamma'})$ is a Hilbert space \cite[pp.~185--187]{Leoni2023}. Regarding the definition of  $H^{1/2}({\Gamma'})$ and its norm, we also reference \cite[Prop.~4.24, p.~192]{Demengel2012}, \cite[Def. 1.3.2.1, pp.~16--17]{Grisvard2011}, \cite[Def. 18.32, 18.36, pp.~613--615]{Leoni2017}, \cite[p.~74]{McLean2000}, \cite[p.~26]{Chouly2023}, \cite[p.~524]{Nezza2012}, \cite[p.~15]{Kikuchi1988}, \cite[Lemma 36.1, p.~196]{Tartar2007}.

\begin{Theorem}{\bf (the trace theorem)} \cite[Th. 9.14, p.~350, Th. 9.39, p.~370]{Leoni2023}, \cite[Th. 3.37, p.~102]{McLean2000}, \cite[Th. 1.3.1, p.~10]{Quarteroni1994} For an open bounded set $\Omega$ with Lipschitz boundary $\Gamma$ there exists a continuous linear operator
\[
{\rm tr}: W^{1,2}(\Omega) \to H^{1/2}(\Gamma)
\]
such that for each $u\in \mathcal{D}\left(\overline{\Omega}\right)$
\[
u(x) = {\rm tr} (u)(x) \qquad \text{ for all }x\in \Omega.
\]
Such an operator is unique. Moreover, ${\rm tr}$ is surjective and has a continuous right inverse:
\[
{\rm rt}: H^{1/2}(\Gamma) \to W^{1,2}(\Omega),
\]
\[
{\rm tr}({\rm rt} (u)) = u \qquad \text{ for all }u\in H^{1/2}(\Omega).
\]
\label{th:trace_th_1}
\end{Theorem}
\begin{Proposition}{\bf} \cite[Th. 18.7, p.~595]{Leoni2017}
\label{prop:W0-ker-tr}
Denote ${\rm Ker}\, {\rm tr}\subset W^{1,2}(\Omega)$ by $W^{1,2}_0(\Omega)$. The space $W^{1,2}_0(\Omega)$ coincides with the closure of $\mathcal{D}(\Omega)$ with respect to the norm of $W^{1,2}(\Omega)$:
\begin{equation}
W^{1,2}_0(\Omega):= {\rm Ker}\, {\rm tr} = \overline{\mathcal{D}(\Omega)}^{\,W^{1,2}(\Omega)}
\label{eq:zero_trace_space_def}
\end{equation}
\end{Proposition}

Throughout the text we will use the notation ${\rm tr}(u)$ for $\mathbb{R}^n$-valued functions, i.e. for $u=(u_k)_{k\in \overline{1,n}} \in[W^{1,2}(\Omega)]^n$ we denote
\begin{equation}
{\rm tr}(u):= ({\rm tr}(u_i))_{i\in \overline{1,n}}\in [H^{1/2}(\Gamma)]^n.
\label{eq:vector-valued-trace-def}
\end{equation}

\subsection{Unbounded adjoint operators}
The framework of {\it unbounded adjoint operators} \cite[Section 2.6]{Brezis2011} will be central for our analysis of the elasticity problem. For convenience, we will include here the main facts to be used.
\begin{Definition}
\label{def:DDC}
{\bf } \cite[p.~43]{Brezis2011} Let $X,Y$ be Banach spaces, $A$ be a linear map from (generally, not closed) subspace  $D(A)$ to $Y$. We write 
\[
A:D(A)\subset X\to Y.
\]
We say that
\begin{enumerate}[{\it i)}]
\item $A$ is {\it densely defined} when $\overline{D(A)}^{X} = X$,
\item $A$ is {\it closed} when for any $\{x_k\}_{k\in \mathbb{N}}\subset D(A),\, x\in X,\, f\in Y$ we have
\[
\begin{cases}
x_k\to x \text{ in }X\\
A x_k \to f \text{ in }Y
\end{cases} \qquad \Longrightarrow\qquad \begin{cases}x\in D(A),\\ f=Ax. \end{cases}
\]
\end{enumerate}
\end{Definition}
\begin{Proposition}{\bf } \cite[Remark 12, p.~43]{Brezis2011}
When $A$ is closed, its kernel ${\rm Ker}\, A$ is a closed subspace of $X$.
\label{prop:closed-operator-closed-kernel}
\end{Proposition}
\begin{Proposition}{\bf } \cite[Def. of $A^*$, Proposition 2.17, pp.~43-45]{Brezis2011} Let $A$ be densely defined. Then there exists a uniquely defined operator
\[
A^*:D(A^*)\subset Y^* \to X^*,
\]
where
\begin{equation}
D(A^*) = \left\{f\in Y^*:\text{there exists }c\geqslant 0: \text{for any }x\in D(A) \text{ we have }|\la f,  A x \ra|\leqslant c\|x\|_X\right\},
\label{eq:D(A)-def}
\end{equation}
such that
\[
\la f, Ax\ra = \la A^* f, x\ra\qquad \text{for any }x\in D(A), f\in D(A^*).
\]
The operator $A^*$ is called {\it the adjoint} of $A$ and $A^*$ is always closed.
\label{prop:adjoint-def}
\end{Proposition}

\begin{Proposition}{\bf } \cite[Corollary 2.18, p.~45]{Brezis2011}
Let $A$ be densely defined and closed. Then 
\begin{enumerate}[{\it i)}]
\item ${\rm Ker}\, A = ({\rm Im}\, A^*)^\perp \xlongequal{def} \left\{x\in X: \text{for any }f\in {\rm Im}\, A^* \text{ we have } \la f, x\ra = 0 \right\},$
\item ${\rm Ker}\, A^* = ({\rm Im}\, A)^\perp \xlongequal{def} \left\{g\in Y^*: \text{for any }y\in {\rm Im}\, A \text{ we have } \la g, y\ra = 0 \right\}$.
\end{enumerate}
\label{prop:prelim-kernel-of-the-adjoint}
\end{Proposition}
\begin{Proposition}{\bf } \cite[Th. 2.19, p.~46]{Brezis2011}
Let $A$ be densely defined and closed. The following properties are equivalent.
\begin{enumerate}[{\it i)}]
\item ${\rm Im}\, A$ is closed,
\item ${\rm Im}\, A^*$ is closed,
\item ${\rm Im}\, A = ({\rm Ker}\, A^*)^\perp \xlongequal{def} \left\{y\in Y: \text{for any }g\in {\rm Ker}\, A^* \text{ we have } \la g, y\ra = 0 \right\}$ ,
\item ${\rm Im}\, A^* = ({\rm Ker}\, A)^\perp \xlongequal{def} \left\{f\in X^*: \text{for any x}\in {\rm Ker}\, A \text{ we have } \la f, x\ra = 0 \right\}$.
\end{enumerate}
\label{prop:prelim-image-of-the-adjoint}
\end{Proposition} 
\begin{Proposition}{\bf } \cite[Th. 2.20, p.~47]{Brezis2011}
Let $A$ be densely defined and closed. The following properties are equivalent.
\begin{enumerate}[{\it i)}]
\item $A$ is surjective,
\item there exists $C\geqslant 0$ such that 
\[
\|g\|_{Y^*} \leqslant C \|A^* g\|_{X^*} \qquad \text{for all }g\in D(A^*).
\]
\item ${\rm Ker}\, A^* = \{0\}$ and ${\rm Im}\, A^*$ is closed.
\end{enumerate}
\label{prop:Brezis-surjectivity}
\end{Proposition}
\begin{Proposition}{\bf }\cite[Th. 3.24, p.~72]{Brezis2011}
Let $X$, $Y$ be reflexive Banach spaces and $A$ be densely defined and closed. Then $A^*$ is densely defined, i.e.
\[
\overline{D(A^*)}^{Y^*}=Y^*
\]
and we have
\[
A^{**}=A.
\]
\label{prop:second-dual}
\end{Proposition}

\subsection{Abstract linear framework}
Given a Hilbert space $\mathcal{H}$ and its closed subspaces $\mathcal{U}, \mathcal{V}$, which are  orthogonal complements of each other in $\mathcal{H}$, we have 
\begin{equation}
\mathcal{U} \oplus \mathcal{V} = \mathcal{H}.
\label{eq:orthogonal_complements_direct_sum}
\end{equation}
Then there exist continuous linear {\it projection operators} $P_{\,\mathcal{U}}, P_{\,\mathcal{V}}$, characterized by
the {\it variational definition of the orthogonal projection}
\begin{equation}
\sigma = P_{\mathcal{V}}g \quad \Longleftrightarrow \quad
\begin{cases} \sigma\in \mathcal{V},\\
\la\tau, \sigma-g\ra=0 & \text{for any }\tau \in \mathcal{V},
\end{cases}
\label{eq:projection_abstract_def}
\end{equation}
and similarly for $P_\mathcal{U}$, see e.g. \cite[Corollary 5.4, p.~134]{Brezis2011}. Also, note the simple yet useful fact that
\begin{equation}
\mathcal{U}\cap\mathcal{V}=\{0\},
\label{eq:abstract_trivial_intersection}
\end{equation}
see e. g. \cite[p.~308]{Rudin1991}.

We will need the following fact about linear operators.
\begin{Proposition}
Let $\mathcal{H}$ be a Hilbert space, split as \eqref{eq:orthogonal_complements_direct_sum} into orthogonal complements $\mathcal{U}, \mathcal{V}$ (possibly, $\mathcal{U}=\mathcal{H}, \mathcal{V}=\{0\}$). Also, let  $\mathcal{X}$ be another Hilbert space, 
\[
A:D(A)\subset \mathcal{H} \to \mathcal{X}
\]
be a densely defined closed linear operator such that
\[
\mathcal{V} = {\rm Ker}\, A,\qquad {\rm Im}\, A \text{ is closed in }\mathcal{X}
\]
(also possibly ${\rm Im}\, A = \mathcal{X}$). Then the restricted operator
\[
A_{\rm r}:D(A_{\rm r})=\left(\mathcal{U}\cap D(A)\right) \subset \mathcal{U} \to {\rm Im}\, A
\]
is also densely defined, closed, bijective and has a continuous inverse
\[
Q:=A_{\rm r}^{-1} \in \mathcal{L}({\rm Im}\, A, \mathcal{U}).
\]
\label{prop:inverse_on_U}
\end{Proposition}
\noindent{\bf Proof.} We refer to \cite[Appendix A]{Beutler1976}. There, in particular, Lemmas A.3 and A.4 imply that $A_{\rm r}$ is closed and densely defined and Lemma A.6 proves the bijection. With the assumption of ${\rm Im}\, A$ being closed we use the implication Lemma A.8 (c) $\Rightarrow$ Lemma A.8 (a) (which follows from the closed graph theorem, see e.g. \cite[Th. 2.9, p. 37]{Brezis2011}) to deduce the boundness of $Q$. $\blacksquare$
\begin{Remark}
The proofs of the facts, which we used for Proposition \ref{prop:inverse_on_U}, remain valid when $\mathcal{X}$ is a Banach space, although they are stated in \cite[Appendix A]{Beutler1976} with $\mathcal{X}$ being a Hilbert space (which is also sufficient for our purposes in the current paper).
\end{Remark}
\begin{Theorem}
In the setting of Proposition \ref{prop:inverse_on_U} further assume that $A$ is surjective. Then for any $g\in \mathcal{H}, f\in \mathcal{X}$ the following problem with unknown $\sigma\in \mathcal{H}$ 
\begin{equation}
\arraycolsep=1pt
\left\{
\begin{array}{rcl}
\sigma&\in& \mathcal{U}+g,\\
A\sigma&=&f
\end{array}
\right.
\label{eq:abstract_linear_problem}
\end{equation}
is equivalent to finding
\begin{equation}
\sigma\in (\mathcal{U}+g) \cap(\mathcal{V}+Qf),
\label{eq:abstract_intersection}
\end{equation}
in which $\mathcal{U}+g=\mathcal{U}+ P_{\,\mathcal{V}}\,g $.
Inclusion \eqref{eq:abstract_intersection} has exactly one solution 
\begin{equation}
\sigma=P_{\,\mathcal{V}}g+Qf,
\label{eq:abstract_linear_solution}
\end{equation} 
which continuously depends on $g$ and $f$ (in the sense of continuity of mapping $\mathcal{H}\times \mathcal{X} \to \mathcal{H}$).
\label{th:abstract_linear_framework_solution}
\end{Theorem}
\noindent{\bf Proof.} By decomposition $g= P_{\,\mathcal{U}}g+ P_{\,\mathcal{V}}g$ we have the equivalence
\[
\sigma\in \mathcal{U}+ g \qquad \Longleftrightarrow\qquad \sigma\in \mathcal{U}+ P_{\,\mathcal{V}}g.
\]
Furthermore,
\[
A\sigma=f \qquad \Longleftrightarrow\qquad \sigma\in {\rm Ker}\,A + Qf = \mathcal{V}+ Qf,
\]
i.e. \eqref{eq:abstract_intersection} is indeed equivalent to \eqref{eq:abstract_linear_problem}.

Notice that \eqref{eq:abstract_intersection} is also equivalent to
\[
\sigma- P_{\,\mathcal{V}}g -Qf \in (\mathcal{U}+P_{\,\mathcal{V}}g) \cap(\mathcal{V}+Qf)-P_{\,\mathcal{V}}g -Qf,
\]
where
\begin{multline*}
(\mathcal{U}+P_{\,\mathcal{V}}g) \cap(\mathcal{V}+Qf)-P_{\,\mathcal{V}}g -Qf  = (\mathcal{U}+P_{\,\mathcal{V}}g-P_{\,\mathcal{V}}g -Qf) \cap(\mathcal{V}+Qf-P_{\,\mathcal{V}}g -Qf)=\\
= (\mathcal{U}-Qf) \cap(\mathcal{V}-P_{\,\mathcal{V}}g)=\mathcal{U}\cap\mathcal{V}=\{0\},
\end{multline*}
where the last two equalities are due to $Qf \in \mathcal{U}, P_{\,\mathcal{V}}\,g\in \mathcal{V}$ and \eqref{eq:abstract_trivial_intersection}, respectively. Therefore, \eqref{eq:abstract_linear_solution} is indeed equivalent to \eqref{eq:abstract_intersection}, and the continuity in \eqref{eq:abstract_linear_solution} follows from continuity of $P_{\,\mathcal{V}}$ and $Q$. $\blacksquare$

\section{Given data, unknown variables and fundamental relations}
\label{sect:datum}
Assume that a body occupies the domain $\Omega$.  We consider two cases:
\begin{itemize}
\item In case of the {\it mixed displacement-traction} boundary conditions
we think of the boundary $\Gamma$ of $\Omega$ as being split into nonempty parts $\Gamma_{\rm D}$ and $\Gamma_{\rm T}$ such that
\begin{gather}
\Gamma_{\rm D} \text{ and }\Gamma_{\rm T}\text{ are open in }\Gamma, \label{eq:boundary-ass-1}
\\
\Gamma_{\rm D} \cap \Gamma_{\rm T} = \varnothing,\label{eq:boundary-ass-2}
\\
\Gamma \setminus (\Gamma_{\rm D}\cup \Gamma_{\rm T}) \text{ is a measure-zero set in }\Gamma. \label{eq:boundary-ass-3}
\end{gather} 
\item In case of the {\it displacement} boundary conditions we have
\[
\Gamma=\Gamma_{\rm D}, \qquad \Gamma_{\rm T}=\varnothing.
\]
\end{itemize}
\begin{Remark} 
\label{remark:Lipschitz-dissection}
In the case of the mixed boundary conditions it is natural to think of a stronger assumption on the ``boundary between boundaries'' than \eqref{eq:boundary-ass-1}--\eqref{eq:boundary-ass-3}, namely, that $\Gamma_{\rm D}$ and $\Gamma_{\rm T}$ satisfy the definition of a Lipschitz domain in $\Gamma$ as the ambient space. But, surprisingly small number of authors pay attention to such a condition. It is briefly mentioned in \cite[Th. 1.3.1, p.~10]{Quarteroni1994}, and the most serious treatment in the literature was provided by McLean \cite[pp.~99,~128,~231]{McLean2000}, who rigorously defined a {\it Lipschitz dissection} of the boundary. 

Somewhat counterintuitively, the abstract problem we consider in this paper, and, in particular, the definitions and the statements of Sections \ref{ssect:displacements-spaces} and \ref{ssect:L-M_spaces} below, do not require the assumption of Lipschitz dissection. In other words, the static elasticity problem with the set $\Gamma \setminus (\Gamma_{\rm D}\cup \Gamma_{\rm T})$ satisfying \eqref{eq:boundary-ass-3}, but still behaving ``badly'' (e.g. the so-called ``Warsaw circle'' \cite[p.~2]{Mardevsic1997}), is a well-defined problem.
\end{Remark}

\subsection{Displacements}
\label{ssect:displacements-spaces}
On the part $\Gamma_{\rm D}$ a Dirichlet boundary condition is prescribed for displacement in the following manner. Assume that 
\[
u_0 \in [W^{1,2}(\Omega)]^n
\]
is given. Consider spaces 
\begin{equation}
\mathcal{D}_{\rm D}: = \{\varphi\in \mathcal{D}\left(\overline{\Omega}\right):{\rm supp}\, \varphi\subset \Omega \cup \Gamma_{\rm T}\},
\label{eq:D-pb-def}
\end{equation}
and its closure in $W^{1,2}(\Omega)$, cf. \eqref{eq:zero_trace_space_def}:
\begin{equation}
W_{\rm D}:=\overline{\mathcal{D}_{\rm D}}^{\,W^{1,2}(\Omega)} = \left\{u\in W^{1,2}(\Omega): {\rm tr}(u)(x)=0 \text{ for a.a. }x\in\Gamma_{\rm D}\right\}.
\label{eq:W-def}
\end{equation}
Observe that
\[
[W^{1,2}_0(\Omega)]^n \subset [W_{\rm D}]^n \subset [W^{1,2}(\Omega)]^n.
\]
Naturally, when $\Gamma_{\rm D}=\Gamma$ we have
\[
[W^{1,2}_0(\Omega)]^n = [W_{\rm D}]^n \subset [W^{1,2}(\Omega)]^n.
\]
\begin{Remark} For a general {\it displacement} $u\in [W^{1,2}(\Omega)]^n$ the norm takes the form
\begin{multline}
\|u\|_{W^{1,2}} =\left(\|u\|_{L^2}^2+\|\nabla u\|_{L^2}^2\right)^{\frac{1}{2}}= \left(\sum_{i=1}^n\|u_i\|_{L^2}^2+\sum_{i,j=1}^n\left\|\frac{\partial u_i}{\partial x_j}\right\|_{L^2}^2\right)^{\frac{1}{2}}=\\ =\left(\sum_{i=1}^n \,\intOm  |u_i|^2 dx+\sum_{i,j=1}^n\,\intOm \left|\frac{\partial u_i}{\partial x_j}\right|^2 dx\right)^{\frac{1}{2}}.
\label{eq:displacement_norms}
\end{multline}
\end{Remark}
\begin{Definition}
\label{def:bc}
We say that an unknown displacement $u\in [W^{1,2}(\Omega)]^n$ satisfies the {\it Dirichlet boundary condition on displacements} associated with given {\it displacement boundary} $\Gamma_{\rm D}$ and {\it displacement load} $u_0\in [W^{1,2}(\Omega)]^n$ when
\[
u-u_0 \in [W_{\rm D}]^n.
\]
Equivalently, we can be given displacement load as $\widetilde{u}_0 \in [H^{1/2}(\Gamma)]^n$ and require that 
\[
{\rm tr}(u) = \widetilde{u}_0 \qquad \text{a.e. on }\Gamma_{\rm D}.
\]
\noindent Respectively, when $\Gamma_{\rm D}=\Gamma$ the Dirichlet boundary condition on displacements is
\[
u-u_0\in [W_0^{1,2}(\Omega)]^n,
\]
or, equivalently,
\[
{\rm tr}(u) = \widetilde{u}_0 \qquad \text{a.e. on }\Gamma.
\]

\end{Definition}

\subsection{Strains}
We will need the space of symmetric matrix-valued fields:
\begin{Definition}
Let 
\begin{multline}
L^2_{\rm sym}(\Omega) =\left\{\tau\in [L^{2}(\Omega)]^{n\times n}: \text{ for all }i,j\in \overline{1,n}\quad \|\tau_{ij}-\tau_{ji}\|_{L^2}=0\right\} = \\ 
= \left\{\tau\in [L^{2}(\Omega)]^{n\times n}: \left\|\tau-\tau^{\top}\right\|_{L^2}=0\right\},
\label{eq:L2sym-def}
\end{multline}
which is a closed subspace of $[L^{2}(\Omega)]^{n\times n}$. Equipped with the inner product and norm of $[L^{2}(\Omega)]^{n\times n}$, space $L^2_{\rm sym}(\Omega)$ is a Hilbert space itself, therefore it is a reflexive Banach space \cite[Proposition 5.1, p.~132]{Brezis2011}.
\label{def:L2sym-def}
\end{Definition}

Using linear operator
\begin{equation}
{\rm E}: [W^{1,2}(\Omega)]^{n} \to L^2_{\rm sym}(\Omega),
\label{eq:E-def1}
\end{equation}
\begin{equation}
{\rm E}:u\mapsto \left(\frac{1}{2} \left(\frac{\partial u_i}{\partial x_j}+ \frac{\partial u_j}{\partial x_i}\right)\right)_{i,j\in\overline{1,n}},
\label{eq:E-def2}
\end{equation}
where the partial derivatives are meant in the weak sense, we define the the {\it displacement-strain relation} between unknowns $u$ and $\varepsilon$ (called {\it strain}) for body $\Omega$ as
\[
\varepsilon = {\rm E}(u)
\]
\begin{Proposition} As defined by \eqref{eq:E-def1}--\eqref{eq:E-def2}, the operator ${\rm E}$ is bounded with operator norm $\|{\rm E}\|_{\rm op}\leqslant 1$.
\label{prop:E-is-bounded}
\end{Proposition}
\noindent{\bf Sketch of the proof.} One directly shows that 
\[
\|{\rm E}(u)\|_{L^2}^2\leqslant \sum_{i,j=1}^n \intOm \left|\frac{\partial u_i}{\partial x_j}\right|^2 dx = \left\|\nabla u\right\|^2_{L^2}.
\] Then \eqref{eq:displacement_norms} implies the statement. $\blacksquare$
\begin{Proposition}{\bf (Korn's first inequality)} \cite[Th. 10.1, p.~298]{McLean2000}, \cite[Theorem 6.3.1, p.~77]{Necas1981}  When $\Omega$ has Lipschitz boundary,
\[
\|\nabla u\|_{L^2}\leqslant\sqrt{2} \|{\rm E}(u)\|_{L^2}\qquad \text{for all }u\in [W^{1,2}_0(\Omega)]^n.
\]
\label{prop:Korn1}
\end{Proposition}
\begin{Proposition}\phantom{V}
\begin{enumerate}[\it i)]
\item {\it (rigid motions)} The kernel of ${\rm E}$ consists of infinitesimal rigid motions (isometries) of the body, i.e.
\[
{\rm Ker\, E} = \{a+Bx: a \in \mathbb{R}^{n}, B\in \mathbb{R}^{n\times n} \text{ such that } B^{\top}=-B\}.
\]
In particular, when $n=3$ we can write
\[
{\rm Ker\, E} = \{a+b\times x: a,b\in \mathbb{R}^{3}\}.
\]
\item When $\Gamma_{\rm D}$ has a positive measure in $\Gamma$, we have 
\[
({\rm Ker\, E}) \cap [W_{\rm D}]^n =\{0\}.
\]
\end{enumerate}
\label{prop:rigid_motions}
\end{Proposition}
\noindent {\bf Proof.} For the proof of {\it i)} we refer to \cite[Lemma 1.1, p.~18]{Temam1985}. To prove {\it ii)} we generalize \cite[Lemma 2.1, p.~91]{Necas1981} from the case $n=3$ to the case $n\geqslant 2$ as follows.

Assume on the contrary that there exists a nonzero function $a+Bx\in [W_{\rm D}]^n$ with $a\in \mathbb{R}^{n}$ and skew-symmetric  $B\in \mathbb{R}^{n\times n}$. If we denote the set
\[
G:=\{x\in \mathbb{R}^n: a+Bx=0\},
\]
then by the assumption we get that $\Gamma_{\rm D}\subset G$. Among all of the possibilities for $a$ and $B$, the only case which allows for $\Gamma_{\rm D}\neq \varnothing$ and $a+Bx\not\equiv0$ is when the  matrix $B$ is nonzero and the equation
\[
Bx=-a
\]
has at least one solution (call it $x^*$). Therefore, 
\begin{equation}
\Gamma_{\rm D}\subset G={\rm Ker}\, B +x^*.
\label{eq:boundary_supset}
\end{equation}
However, $B$ is a nonzero skew-symmetric matrix, thus ${\rm rank}\, B\geqslant 2$. Indeed,
$B$ has an element $b_{ij}\neq 0$ with $i\neq j$ and, correspondingly, $b_{ji}=-b_{ij}\neq 0$. By renumbering the components we can bring the matrix to a form with $i=1, j=2$, i.e. without the loss of generality the upper left $2\times 2$ submatrix of $B$ is of the form
\[
\begin{pmatrix}
0& b\\
-b &0
\end{pmatrix}
\]
with some $b\neq 0$. Hence, $B$ has at least two linearly independent rows.

By the rank-nullity theorem (see e.g. \cite[Th. 2.49, p.~114]{OlverLinearAlgebra2018}) 
\[
{\rm dim\, Ker}\, B = n-{\rm rank}\, B \leqslant n-2.
\]
But this is impossible, because $\Gamma_{\rm D}$ in \eqref{eq:boundary_supset} has a positive measure in $\Gamma$, which is of dimension $n-1$. The contradiction means that $a+Bx\equiv 0$. $\blacksquare$

\begin{Proposition} If $\Gamma_{\rm D}$ has a positive measure in $\Gamma$, then there exists $c>0$ such that
\begin{equation}
\|u\|_{W^{1,2}}\leqslant c\, \|{\rm E}(u)\|_{L^2} \qquad \text{for all }u\in [W_{\rm D}]^n.
\label{eq:fromKorn2}
\end{equation}
\end{Proposition}
\noindent {\bf Proof.} We follow \cite[Section 3.3.1, pp.~78-81]{Chouly2023}. The estimate \eqref{eq:fromKorn2} is a direct consequence of {\it Korn's second inequality} (see e.g. \cite[Th. 10.2, p.~299]{McLean2000}, \cite[Th. 5.13, p.~106]{Kikuchi1988}, \cite[(3.7), p.~79, Th. 3.4, p.~84]{Necas1981}) by Peetre--Tartar lemma (see \cite[Lemma 3.5, p.~79]{Chouly2023}, \cite[Lemma A.38, p.~469]{Ern2004}). Peetre--Tartar lemma requires an appropriate choice of a bounded injective operator $A$ and a compact operator $T$. 

The operator $A$ is chosen as the restriction ${\rm E}|_{[W_{\rm D}]^n}$, which is bounded by construction, see Proposition \ref{prop:E-is-bounded}. The kernel of ${\rm E}$ intersected with $[W_{\rm D}]^n$ is trivial when $\Gamma_{\rm D}$ is of positive measure (see Proposition \ref{prop:rigid_motions}), hence ${\rm E}|_{[W_{\rm D}]^n}$ is injective.

In turn, $T$ is chosen as the embedding of $[W_{\rm D}]^n$ into $[L^2(\Omega)]^n$, which is compact by Rellich's theorem (see e.g. \cite[Th. 2.3, p.~75]{Necas1981}, \cite[Th. 3.27, p.~87]{McLean2000}). $\blacksquare$

\noindent For alternative presentations of the proof of \eqref{eq:fromKorn2} also see \cite[Def. 3.1, p.~79, Th. 3.5, p.~85]{Necas1981}, \cite[Lemma 6.2, p.~115]{Kikuchi1988}.

\subsection{Stresses}
\label{sssect:stresses}
Let for a.a. $x\in \Omega$ we are given a linear map (called the {\it elasticity tensor})
\begin{equation}
{\bf C}(x): \mathbb{S} \to \mathbb{S},
\label{eq:C-def1}
\end{equation}
where $\mathbb{S}$ is the space of symmetric $n\times n$ matrices:
\[
\mathbb{S} = \left\{(\tau_{ij})_{i,j\in \overline{1,n}}\in \mathbb{R}^{n\times n}: \tau_{ij}=\tau_{ji} \text{ for all }i,j\in \overline{1,n}\right\}.
\]
We assume that ${\bf C}$ satisfies the following conditions:
\begin{enumerate}[\it i)]
\item {\it ($L^{\infty}$ -regularity and uniform coercivity)} There exist $c_0, C_0>0$ such that 
\begin{equation}
c_0\, \tau:\tau \leqslant \tau: ({\bf C}(x)\,  \tau) \leqslant C_0\, \tau:\tau \qquad \text{for all } \tau \in \mathbb{S}\text{ and a.a. }x\in \Omega.
\label{eq:bounded_and_coercive}
\end{equation}
\item{\it (symmetry)}
\begin{equation}
\tau_1 :  ({\bf C}(x)\,  \tau_2) =  ({\bf C}(x)\,  \tau_1):\tau_2   \qquad \text{for all }\tau_1, \tau_2\in \mathbb{S}\text{ and a.a. }x\in \Omega.
\label{eq:symmetric}
\end{equation}
\end{enumerate}
Then {\it Hooke's law}, connecting unknowns $\varepsilon$ and $\sigma$ (the latter is called {\it stress}) at point $x\in \Omega$ takes the form
\begin{equation}
\sigma(x) = {\bf C}(x)\, \varepsilon(x).
\label{eq:C-def2}
\end{equation}
We interpret this as relation $\sigma = {\bf C}\, \varepsilon$ in $L^{2}_{\rm sym}(\Omega)$, where
\begin{equation}
{\bf C}: L^{2}_{\rm sym}(\Omega) \to L^{2}_{\rm sym}(\Omega),
\label{eq:C-def3}
\end{equation}
which is well-defined due to the right estimate of \eqref{eq:bounded_and_coercive}. Moreover, from \eqref{eq:bounded_and_coercive}--\eqref{eq:symmetric} one can deduce that the inverse $\left({\bf C}(x)\right)^{-1}$ exists for a.a. $x\in \Omega$ and satisfies 
\begin{gather}
C^{-1}_0\, \tau:\tau \leqslant \tau: ({\bf C}^{-1}(x)\,  \tau) \leqslant c^{-1}_0\, \tau:\tau \qquad \text{for all } \tau \in \mathbb{S}\text{ and a.a. }x\in \Omega, \label{eq:bounded_and_coercive_C_inv}\\[1mm]
\tau_1 :  ({\bf C}^{-1}(x)\,  \tau_2) =  ({\bf C}^{-1}(x)\,  \tau_1):\tau_2   \qquad \text{for all }\tau_1, \tau_2\in \mathbb{S}\text{ and a.a. }x\in \Omega.\label{eq:symmetric_C_inv}
\end{gather}
Hence, we can also write $\varepsilon={\bf  C}^{-1}\sigma$ as an equality in $L^{2}_{\rm sym}(\Omega)$ with
\[
{\bf C}^{-1}:L^{2}_{\rm sym}(\Omega) \to L^{2}_{\rm sym}(\Omega).
\]
\subsection{Force loads and the equation of equilibrium}
\label{ssect:eq-equil-def}
Let us have prescribed {\it body force load} over $\Omega$ and {\it traction force load} over $\Gamma$ as given functions, respectively,
\[
F\in [L^{2}(\Omega)]^n, \qquad T\in [L^2(\Gamma)]^n.
\]
As we will discuss below, one can consider $T$ defined only over $\Gamma_{\rm T}$ or discard it completely when $\Gamma_{\rm T}=\varnothing$, i.e. it is enough to consider
\begin{equation}
T\in [L^2(\Gamma_{\rm T})]^n.
\label{eq:traction_force_datum}
\end{equation}

The {static problem of elasticity} consists of finding the values of the unknown variables $u, \varepsilon, \sigma$, corresponding to an {\it equilibrium} state of the body. As a criterion of the equilibrium the {\it principle of virtual work} can be used (for continuum mechanics see e.g. \cite[Section 2.2.2, pp.~63--65]{Altenbach2018}, \cite[Appendix 1, pp.~589--591]{Anandarajah2011}, \cite[Remark 5.1.1, p.~61]{Necas1981}, cf. \cite[Section 1.4, pp.~16--17]{Goldstein2002} for systems of particles). In particular, for the model we described so far it takes the following form.

\begin{Definition} 
\label{def:virtual_work_equilibrium}
We say that 
\[
(u, \varepsilon, \sigma)\in [W^{1,2}(\Omega)]^3\times L^2_{\rm sym}(\Omega)\times L^2_{\rm sym}(\Omega)
\]
 corresponds to a state of equilibrium when for any {\it virtual displacement} $u+\delta u$ satisfying the Dirichlet boundary condition, the {\it virtual work} that would be done by the {\it internal forces} over the transition from the configuration with the displacement $u$ to the configuration with the displacement $u+\delta u$ equals the virtual work that would be done by the {\it external forces} over the same transition, see the references above. Specifically, this means that for any $u+\delta u$ satisfying the condition of Definition \ref{def:bc} we must have
\[
{\rm IVW}(\delta u) = {\rm EVW} (\delta u),
\]
where the {\it internal virtual work} (IVW) and {\it external virtual work} (EVW) are, respectively,
\begin{align*}
{\rm IVW}(\delta u)&= \intOm \sigma:({\rm E}(u+\delta u) - {\rm E}(u))\, dx,
\\
{\rm EVW}(\delta u) &= \intOm F\cdot (u+\delta u - u)\,dx + \int\limits_{\Gamma} T\cdot {\rm tr}\,(u+\delta u-u) \,dS.
\end{align*}
\end{Definition} 

Due to the fact that both $u$ and $u+\delta u$ must satisfy the conditions of Definition \ref{def:bc} and because of linearity of ${\rm E}$, Definition \ref{def:virtual_work_equilibrium} yields the following {\it equation of equilibrium}
\[
\intOm \sigma: {\rm E}(\delta u)\, dx = \intOm F\cdot \delta u \,dx + \int\limits_{\Gamma_{\rm T}} T\cdot {\rm tr} (\delta u) \,dS \qquad \text {for any }\delta u\in [W_{\rm D}]^n,
\]
or, when $\Gamma_{\rm T}=\varnothing$:
\[
\intOm \sigma: {\rm E}(\delta u)\, dx = \intOm F\cdot \delta u \,dx  \qquad \text {for any }\delta u\in [W^{1,2}_0(\Omega)]^n.
\]
Note that the part of the integral with $T$ over $\Gamma_{\rm D}$ vanished (due to ${\rm tr}(\delta u)=0$ on $\Gamma_{\rm D}$) and $T$ does not appear in other equations. Thus the values of $T$ on $\Gamma_{\rm D}$ have no effect on the model, and the datum can be given as in \eqref{eq:traction_force_datum} without the loss of generality.

\subsection{The complete statements of the problems}
For convenience we combine the formulations of the problems in the following definitions.
\begin{Definition}
\label{def:static-problem-of-linear-elasticity-mixed}
We say that the unknowns of 
\begin{align}
\text{displacement field } & u \in [W^{1,2}(\Omega)]^n, \label{eq:displacement-variable-from}\\ 
\text{strain field } & \varepsilon \in L^{2}_{\rm sym}(\Omega),\label{eq:strain-variable-from}\\ 
\text{stress field } & \sigma \in  L^{2}_{\rm sym}(\Omega) \label{eq:stress-variable-from}
\end{align}
solve the {\it static problem of linear elasticity with mixed displacement-traction boundary conditions} with prescribed
\begin{align*}
\text{displacement }& u_0 \in [W^{1,2}(\Omega)]^n \text { on the part of the boundary }\Gamma_{\rm D},\\ 
\text{traction force load } & T\in [L^2(\Gamma_{\rm T})]^n  \text{ on the part of the boundary } \Gamma_{\rm T},\\
\text{body force load } & F\in [L^{2}(\Omega)]^n 
\end{align*}
when the following holds:
\begin{align} \text{Dirichlet boundary condition: }& u-u_0\in [W_{\rm D}]^n, \label{eq:le-1} \tag{LE1} \\
\text{displacement-strain relation: }& \varepsilon = {\rm E}(u), \label{eq:le-2} \tag{LE2} \\
\text{Hooke's law: }& \sigma = {\bf C}\, \varepsilon, \label{eq:le-3} \tag{LE3} \\
\text{equation of equilibrium: }& \intOm \sigma: {\rm E}(\delta u)\, dx = \intOm F\cdot \delta u \,dx + \int\limits_{\Gamma_{\rm T}} T\cdot {\rm tr} (\delta u) \,dS \quad \text {for any }\delta u\in [W_{\rm D}]^n, \tag{LE4} \label{eq:le-4}
\end{align}
where $W_{\rm D}$ is defined by \eqref{eq:W-def}, ${\rm E}$ by \eqref{eq:E-def1}--\eqref{eq:E-def2} and ${\bf C}$ by \eqref{eq:C-def1}--\eqref{eq:C-def3}.
\end{Definition}
\begin{Definition}
\label{def:static-problem-of-linear-elasticity-displacement}
In turn, the triple $(u, \varepsilon, \sigma)$ as in \eqref{eq:displacement-variable-from}--\eqref{eq:stress-variable-from} solves the {\it static problem of linear elasticity with the displacement boundary conditions} with prescribed 
\begin{align*}
\text{displacement }& u_0 \in [W^{1,2}(\Omega)]^n \text { on the boundary }\Gamma,\\ 
\text{body force load } & F\in [L^{2}(\Omega)]^n 
\end{align*}
when \eqref{eq:le-2}--\eqref{eq:le-3} hold together with
 \begin{align} \text{Dirichlet boundary condition: }& u-u_0\in [W^{1,2}_{0}(\Omega)]^n, \label{eq:le-1'} \tag{LE1$'$} \\
\text{equation of equilibrium: }& \intOm \sigma: {\rm E}(\delta u)\, dx = \intOm F\cdot \delta u \,dx \quad \text {for any }\delta u\in [W^{1,2}_0(\Omega)]^n. \tag{LE4$'$} \label{eq:le-4'}
\end{align}
\end{Definition}
In this text we pursue the solutions of the following {\it stress problems}.
\begin{Definition} 
\label{def:reduced-problem-of-linear-elasticity}
We say that $\sigma \in  L^{2}_{\rm sym}(\Omega)$ solves the {\it static stress problem of linear elasticity with the mixed boundary conditions} (resp.~{\it with the displacement boundary conditions}) when there exist $u, \varepsilon$ as in \eqref{eq:displacement-variable-from}--\eqref{eq:strain-variable-from} so that the triple $(u, \varepsilon, \sigma)$ satisfies the conditions of Definition \ref{def:static-problem-of-linear-elasticity-mixed} (resp.~Definition \ref{def:static-problem-of-linear-elasticity-displacement}) with corresponding given data.
\end{Definition}

\section{Analysis of the static stress problem of linear elasticity with the displacement boundary conditions}
\label{sect:problem1}
\subsection{Preliminaries: weak divergence and the associated spaces}
\label{ssect:weak_divergence}
In order to work with the equations \eqref{eq:le-4}, \eqref{eq:le-4'} it is helpful to introduce the Hilbert spaces of functions with weak divergence. We reference \cite[(1.7)--(1.8), p.~38]{Krizek1983} and \cite[(1.43),(1.45), p.~14]{Temam1985} for symmetric matrix-valued functions, \cite[p.~40]{Gatica2014} for matrix-valued functions, and \cite[Section 1.2.1, p.~13]{Girault1979}, \cite[Section 3.6, p.~40]{Krizek1996}, \cite[p.~27]{Chouly2023}, \cite[p.~15]{Gatica2014}, \cite[pp.~108--109]{Sayas2019}, \cite[p.~7, pp.~49--51]{Boffi2013} for $\mathbb{R}^n$-valued functions.

\begin{Definition}
If functions $\tau\in [L^2(\Omega)]^{n\times n}$ and $f\in [L^2(\Omega)]^n$ satisfy
\begin{equation}
\intOm \tau:\nabla u\, dx = -\intOm f\cdot u\, dx \qquad \text{for all } u\in [\mathcal{D}(\Omega)]^n 
\quad(\text{equivalently, for all } u\in [W^{1,2}_0(\Omega)]^n),
\label{eq:div-def}
\end{equation}
we say that $\tau\in H({\rm div}, \Omega)$ and $f$ is the weak divergence of $\tau$, denoted by ${\rm div}\, \tau$. Additionally, let
\begin{equation}
H_{\rm sym}({\rm div}, \Omega) = H({\rm div}, \Omega) \cap L^2_{\rm sym}(\Omega),
\label{eq:H_div_sym-def}
\end{equation}
\label{def:H_div_sym-def}
\end{Definition}

\begin{Proposition}
Equipped with the following inner product and norm, the space $H({\rm div}, \Omega)$ is a Hilbert space:
\begin{equation}
\la \tau_1, \tau_2 \ra_{H({\rm div}, \Omega)}=\la \tau_1 ,\tau_2 \ra_{L^2} + \left\la {\rm div}\, \tau_1, {\rm div} \,\tau_2 \right\ra_{L^2} = \intOm\left( \tau_1 :\tau_2 +{\rm div}\, \tau_1 \cdot  {\rm div} \,\tau_2\right) \,dx,
\label{eq:H-div-inner-product}
\end{equation}
\begin{equation}
\left\| \tau\right\|_{H({\rm div}, \Omega)} = \left(\la\tau,\tau \ra_{H({\rm div}, \Omega)}\right)^{\frac{1}{2}}.
\label{eq:H-div-norm}
\end{equation}
The space $H_{\rm sym}({\rm div}, \Omega)$ is also a Hilbert space with inner product \eqref{eq:H-div-inner-product} and norm \eqref{eq:H-div-norm}.
\label{prop:H-div-Hilbert}
\end{Proposition}
\noindent {\bf Proof.}
Clearly, \eqref{eq:H-div-inner-product} is symmetric, bilinear and positive definite,  so we only need to show the completeness of $H({\rm div}, \Omega)$. Assume that $\{\tau_k\}_{k\in \mathbb{N}}\subset H({\rm div}, \Omega)$ is a Cauchy sequence in norm \eqref{eq:H-div-norm}. Then there are $\tau^*\in [L^2(\Omega)]^{n\times n}, f^*\in [L^2(\Omega)]^n$ such that
\[
\|\tau_k-\tau^*\|_{L^2}\to 0, \qquad \|{\rm div}\, \tau_k - f^*\|_{L^2} \to 0.
\]
Observe that by Definition \ref{def:H_div_sym-def} we have for any $u\in [\mathcal{D}(\Omega)]^n$
\begin{align*}
\left|\intOm\tau^*:\nabla u\, dx + \intOm f^*\cdot u\, dx \right| &= \left|\intOm(\tau^*- \tau_k):\nabla u\, dx +\intOm \tau_k:\nabla u\, dx+ \intOm f^*\cdot u\, dx \right| =\\
&= \left|\intOm(\tau^*- \tau_k):\nabla u\, dx + \intOm (f^*-{\rm div}\, \tau_k)\cdot u\, dx \right|\leqslant \\[1mm]
&\leqslant\left\|\tau^*-\tau_k\right\|_{L^2}\, \|\nabla u\|_{L^2} + \left\|f^*-{\rm div}\, \tau_k\right\|_{L^2}\, \|u\|_{L^2}\to 0,
\end{align*}
and therefore $H({\rm div}, \Omega)$ is a Hilbert space. The statement on $H_{\rm sym}({\rm div}, \Omega)$ directly follows from the fact that it is closed in $H({\rm div}, \Omega)$. $\blacksquare$

\begin{Proposition}
\label{prop:E-grad-symmetry-equality}
For any $u\in [W^{1,2}(\Omega)]^n$ and $\tau\in  L^2_{\rm sym}(\Omega)$ we have
\begin{equation}
\intOm \frac{1}{2}\left(\nabla u+ (\nabla u)^{\top}\right):\tau \, dx =\intOm \tau:\nabla u\, dx. 
\label{eq:E-grad-symmetry-equality}
\end{equation}
\end{Proposition}
\noindent {\bf Proof.} Observe that 
\begin{equation}
\begin{aligned}
\left|\intOm \frac{1}{2} (\nabla u) : \tau\, dx -\intOm \frac{1}{2} (\nabla u)^{\top} : \tau \, dx  \right|&= \left|\intOm \frac{1}{2} (\nabla u): \tau\, dx -\intOm \frac{1}{2} (\nabla u) : \tau^{\top}\, dx  \right|=\\
&=\frac{1}{2}\left|\intOm (\nabla u): (\tau-\tau^{\top})\, dx\right|\leqslant \frac{1}{2} \|\nabla u\|_{L^2}\,\|\tau-\tau^{\top}\|_{L^2} = 0
\end{aligned}
\label{eq:E-grad-symmetry-equality-proof}
\end{equation}
due to $\tau$ being from \eqref{eq:L2sym-def}. The statement of the proposition follows from the summation of the left-hand sides of \eqref{eq:E-grad-symmetry-equality} and \eqref{eq:E-grad-symmetry-equality-proof} (without the absolute value). $\blacksquare$
\subsection{The pair of adjoint operators for the problem with the displacement boundary conditions}
\label{ssect:aux_operators}
Consider the operators
\begin{gather}
{\bf E_0}: D({\bf E_0}) = [W_0^{1,2}(\Omega)]^n \subset [L^2(\Omega)]^n \to L^2_{\rm sym}(\Omega), \label{eq:E0-def1}\\
{\bf E_0}:u \mapsto \frac{1}{2}\left(\nabla u+ (\nabla u)^{\top}\right),
\label{eq:E0-def2}
\end{gather}
\begin{gather}
{\bf D_0}:D({\bf D_0}) = H_{\rm sym}({\rm div}, \Omega) \subset  L^2_{\rm sym}(\Omega)\to [L^2(\Omega)]^n \label{eq:D0-def1},
\\
{\bf D_0}: \sigma\mapsto - {\rm div}\, \sigma.
\label{eq:D0-def2}
\end{gather}
\begin{Theorem}
\label{th:auxiliary_operators_are_mutually_adjoint}
Both ${\bf E_0}$ and ${\bf D_0}$ are closed densely defined unbounded operators (see Definition \ref{def:DDC}),  and they are mutually adjoint operators, i.e. in the notation of Proposition \ref{prop:adjoint-def} we have
\begin{equation}
{\bf E_0}^* = {\bf D_0} \quad \text{ and }\quad {\bf D_0}^* = {\bf E_0}.
\label{eq:E0-D0-mutually-adjoint}
\end{equation}
In particular, the density of the domain of ${\bf D_0}$ (Definition \ref{def:DDC} {\it i}) means that
\begin{equation}
\overline{H_{\rm sym}({\rm div}, \Omega)}^{[L^2(\Omega)]^{n\times n}}=L^2_{\rm sym}(\Omega).
\label{eq:density_H_sym_in_L2}
\end{equation}
\end{Theorem}
\noindent {\bf Proof.} \newline
\underline{Claim: ${\bf E_0}$ is densely defined.}
It follows from the fact that $[D(\Omega)]^n \subset [W_0^{1,2}(\Omega)]^n$ is dense in $[L^2(\Omega)]^n$, see e.g. \cite[Corollary 4.23, p.~109]{Brezis2011}. 
\newline \underline{Claim: ${\bf E_0}$ is closed.} Assume that 
\[
\{u_k\}_{k\in \mathbb{N}}\subset [W_0^{1,2}(\Omega)]^n, \quad u^*\in [L^2(\Omega)]^n, \quad f \in L^2_{\rm sym}(\Omega)
\] 
are such that
\begin{equation}
\|u_k-u^*\|_{L^2}\to 0,
\label{eq:E-cl1}
\end{equation}
\begin{equation}
\left\|{\bf E_0}\, u_k - f\right\|_{L^2} \to 0.
\label{eq:E-cl2}
\end{equation}
Notice that $\nabla u_{k}\in [L^2(\Omega)]^{n\times n}$ are well-defined for each $k\in \mathbb{N}$. Moreover, by Korn's first inequality (Proposition \ref{prop:Korn1})
\[
\|\nabla u_k\|_{L^2}\leqslant \sqrt{2}\left\|{\bf E_0}\, u_k\right\|_{L^2}.
\]
From \eqref{eq:E-cl2} we deduce the boundness, i.e. that there exists $C_1\geqslant 0$ such that for any $k\in \mathbb{N}$
\[
\left\|{\bf E_0}\, u_k - f \right\|_{L^2} \leqslant C_1,
\]
since the left-hand side is a Cauchy sequence. Thus
\[
\|\nabla u_k\|_{L^2}\leqslant \sqrt{2}\left\|{\bf E_0} \,u_k - f \right\|_{L^2} +\sqrt{2}\|f\|_{L^2} \leqslant \sqrt{2}C_1 + \sqrt{2}\|f\|_{L^2} =:C_2.
\]
By Kakutani's theorem \cite[Th. 3.17, p.~67]{Brezis2011} and the reflexivity of $[L^2(\Omega)]^{n\times n}$ there exists a subsequence $\{u_{k_l}\}_{l\in \mathbb{N}}\subset \{u_k\}_{k\in\mathbb{N}}$ and the weak limit $g\in [L^2(\Omega)]^{n\times n}$ of gradients:
\begin{equation}
\left|\intOm \left(\nabla u_{k_l} - g\right):\tau\, dx\right| \to 0 \text{ as }l\to \infty \qquad\text{ for any }\tau\in [L^2(\Omega)]^{n\times n}.
\label{eq:weak-convergence-subsequence}
\end{equation}
Hence,
\begin{equation}
\left|\intOm \left(\frac{\partial u_{ik_l}}{\partial x_j} - g_{ij}\right) \varphi\, dx\right| \to 0 \text{ as }l\to \infty \qquad\text{ for any }\varphi \in \mathcal{D}(\Omega) \text{ and any } \, i, j\in \overline{1,n},
\label{eq:weak-convergence-components}
\end{equation}
where $u_{ik}$ are the components of $u_k$ and $g_{ij}$ are the components of $g$.
Notice that 
\begin{align*}
\left|\intOm \left(g_{ij} - \frac{\partial u_{ik_l}}{\partial x_j}\right)\varphi \,dx\right| &=  \left|\intOm g_{ij}\varphi\, dx + \intOm u_{ik_l} \frac{\partial \varphi}{\partial x_j} \,dx\right| =\\
&= \left|\intOm g_{ij}\varphi\, dx + \intOm u_i^* \frac{\partial \varphi}{\partial x_j}\, dx - \intOm \left(u_i^* - u_{ik_l}\right) \frac{\partial \varphi}{\partial x_j} \,dx\right| \geqslant \\
&\geqslant
\left|\intOm g_{ij}\varphi\, dx + \intOm u_i^* \frac{\partial \varphi}{\partial x_j}\, dx\right| - \left|\intOm \left(u_i^* - u_{ik_l}\right) \frac{\partial \varphi}{\partial x_j} \,dx\right|,
\end{align*}
where $u^*_i$ is the $i$-th component of $u^*$. Thus,
\begin{align*}
\left|\intOm g_{ij}\varphi\, dx + \intOm u_i^* \frac{\partial \varphi}{\partial x_j}\, dx\right| &\leqslant \left|\intOm \left(g_{ij} - \frac{\partial u_{ik_l}}{\partial x_j}\right)\varphi \,dx\right| + \left|\intOm \left(u_i^* - u_{ik_l}\right) \frac{\partial \varphi}{\partial x_j} \,dx\right|\leqslant\\
&\leqslant \left|\intOm \left(g_{ij} - \frac{\partial u_{ik_l}}{\partial x_j}\right)\varphi \,dx\right| + \left\|u^*_i - u_{ik_l}\right\|_{L^2} \left\| \frac{\partial \varphi}{\partial x_j}\right\|_{L^2} \to 0,
\end{align*}
where the zero limit is due to \eqref{eq:weak-convergence-components} and \eqref{eq:E-cl1}.
We have obtained
\[
\intOm g_{ij}\varphi\, dx + \intOm u_i^* \frac{\partial \varphi}{\partial x_j}\, dx = 0 \qquad \text{ for any }\varphi \in \mathcal{D}(\Omega)
\]
with $g_{ij}\in L^2(\Omega)$ for each $i,j\in \overline{1,n}$. By the definition, $g_{ij}$ is the weak derivative $\frac{\partial u^*_i}{\partial x_j}$ and 
\[
u^*\in [W^{1,2}(\Omega)]^n
\]
with
\begin{equation}
\nabla u^* = g.
\label{eq:proof-weak-gradient-of-the-limit}
\end{equation}
Additionally, it follows that $g$ is independent of the choice of the subsequence $u_{k_l}$ due to the uniqueness of the weak derivative (see e.g. \cite[p.~257]{Evans2010}).

Observe that 
\begin{equation}
u_{k_l} \to u^* \qquad \text{weakly in }[W^{1,2}(\Omega)]^n.
\label{eq:weak-convergence-in-W12}
\end{equation}
Indeed, for any $v\in [W^{1,2}(\Omega)]^n$
\[
\left| \intOm \left(u_{k_l}-u^*\right)\cdot v\, dx +  \intOm \left(\nabla u_{k_l}-g\right): \nabla v\,dx  \right|\leqslant \left\|u_{k_l} - u^*\right\|_{L^2} \|v\|_{L^2} + \left|\intOm \left(\nabla u_{k_l}-g\right): \nabla v\,dx\right| \to 0,
\]
where the zero limit is due to \eqref{eq:E-cl1} and \eqref{eq:weak-convergence-subsequence}. 

Since $[W^{1,2}_0(\Omega)]^n$ is a convex set, closed in $[W^{1,2}(\Omega)]^n$, it is also weakly closed, see \cite[Th. 3.7, p.~60]{Brezis2011}. Recall that we took $u_k\in [W^{1,2}_0(\Omega)]^n$, and hence \eqref{eq:weak-convergence-in-W12} yields
\[
u^*\in [W^{1,2}_0(\Omega)]^n = D({\bf E_0}).
\]
At last, observe that 
\[
{\bf E_0}\, u_{k_l} \to {\bf E_0} \,u^* \qquad \text{weakly in } L^2_{\rm sym}(\Omega).  
\]
Indeed, for any $\tau\in L^2_{\rm sym}(\Omega)$
\begin{align*}
\left| \intOm \left({\bf E_0}\, u_{k_l} - {\bf E_0}\, u^*\right):\tau\,  dx\right| &= \left| \intOm {\bf E_0}(u_{k_l} - u^*):\tau\,  dx\right| = \\
&=\left| \intOm \nabla(u_{k_l} - u^*):\tau\, dx\right| = \left| \intOm \left(\nabla u_{k_l} - g\right):\tau\, dx\right| \to 0,
\end{align*}
where the second equality is due to Proposition \ref{prop:E-grad-symmetry-equality}, the last equality is due to \eqref{eq:proof-weak-gradient-of-the-limit}, and the zero limit is due to \eqref{eq:weak-convergence-subsequence}.

On the other hand, norm-limit \eqref{eq:E-cl2} implies that 
\[
{\bf E_0}\, u_{k_l} \to f\qquad \text{weakly in } L^2_{\rm sym}(\Omega). 
\]
as well. By the uniqueness of the weak limit (see e. g. \cite[p.~232]{BachmanNarici1966})
\[
f= {\bf E_0}\, u^*.
\]
We have shown that ${\bf E_0}$ is closed.
\newline \underline{Claim: $D\left({\bf E_0}^*\right) \subset D\left({\bf D_0}\right)$.} As defined by \eqref{eq:D(A)-def},
\begin{multline}
D\left({\bf E_0}^*\right) =\\ =\left\{\sigma\in L^2_{\rm sym}(\Omega): \text{ there exists }c\geqslant 0: \text{ for all }u\in [W_0^{1,2}(\Omega)]^n \vphantom{\left|\intOm\right|} \right. \quad \left.
 \left|\intOm \sigma: {\bf E_0}\, u \, dx \right|\leqslant  c \|u\|_{L^2}  \right\}.
\label{eq:D(E*)-def}
\end{multline}
In other words, each $\sigma\in L^2_{\rm sym}(\Omega)$ defines a linear functional $\widetilde{\sigma}$
\[
\widetilde{\sigma}: D\left(\widetilde{\sigma}\right)=[W_0^{1,2}(\Omega)]^n\subset [L^2(\Omega)]^n \to \mathbb{R},
\]
\begin{equation*}
\widetilde{\sigma}: u\mapsto \intOm \sigma: {\bf E_0}\, u\, dx,
\end{equation*}
for which $\sigma\in D({\bf E_0}^*)$ means boundness with respect to the $[L^2(\Omega)]^n$-norm in $D\left(\widetilde{\sigma}\right)$. The functional $\widetilde{\sigma}$ can be uniquely extended to be defined and bounded on the entire $[L^2(\Omega)]^n$, see \cite[p.~44]{Brezis2011}, and we denote the extension by $\widetilde{w}\in \left([L^2(\Omega)]^n\right)^*$. In turn, there exists a Riesz representation $w\in[L^2(\Omega)]^n$ of $w^*$ such that
\[
\la \widetilde{w}, u\ra = \intOm w \cdot u \, dx \qquad \text{ for any }u\in[L^2(\Omega)]^n.
\]
As $\sigma\in L^2_{\rm sym}(\Omega)$, from the definition of ${\bf E_0}$ and Proposition \ref{prop:E-grad-symmetry-equality} we conclude that
\[
\intOm \sigma: \nabla u\,dx = \intOm \sigma: {\bf E_0}\, u\,dx = \la\widetilde{\sigma},u\ra = \la\widetilde{w},u\ra = \intOm w\cdot u\, dx \qquad \text{ for all }u\in [W_0^{1,2}(\Omega)]^n.
\]
Observe that the pair $(\sigma, -w)$ satisfies Definition \ref{def:H_div_sym-def} in the sense that $\sigma\in H({\rm div}, \Omega)$ and $-w={\rm div}\, \sigma$. Recall again that $\sigma\in L^2_{\rm sym}(\Omega)$, and therefore $\sigma\in H_{\rm sym}({\rm div}, \Omega)$ accordingly to how it is defined by \eqref{eq:H_div_sym-def}. We have shown that
\[
D({\bf E_0}^*) \subset D({\bf D_0}).
\]

\noindent\underline{Finishing the proof.}
Vice versa, take any $\sigma\in H_{\rm sym}({\rm div},\Omega)$ and observe that for any \newline $u\in [W_0^{1,2}(\Omega)]^n$
\begin{equation}
\intOm \sigma: {\bf E_0} \, u\,dx  = \intOm \sigma: \nabla u\,dx = \intOm(-{\rm div}\, \sigma)\cdot u\, dx = \intOm({\bf D_0}\, \sigma)\cdot u\, dx
\label{eq:E0-Div_sym-duality}
\end{equation}
thus
\[
\left| \intOm \sigma: {\bf E_0}\, u\,dx \right| =\left| \intOm(-{\rm div}\, \sigma)\cdot u\, dx\right|\leqslant \|{\rm div}\, \sigma\|_{L^2}\|u\|_{L^2} \qquad \text{ for any }   u\in [W_0^{1,2}(\Omega)]^n,
\]
i.e. $\sigma\in D({\bf E_0}^*)$ according to \eqref{eq:D(E*)-def}. We have established that
\[
D({\bf E_0}^*) = H_{\rm sym}({\rm div}, \Omega) =  D({\bf D_0}).
\]
Since \eqref{eq:E0-Div_sym-duality} holds  for any $u\in [W_0^{1,2}(\Omega)]^n$ and any $\sigma\in H_{\rm sym}({\rm div}, \Omega)$, we deduce that
\[
{\bf E_0}^* = {\bf D_0}
\]
and ${\bf D_0}$ is a closed operator (see Proposition \ref{prop:adjoint-def}). Moreover, since both $L^2_{\rm sym}(\Omega)$ and $[L^2(\Omega)]^n$ are reflexive, we have that $ {\bf D_0}$ is densely defined (i.e. \eqref{eq:density_H_sym_in_L2} holds) and 
\[
{\bf D_0}^*={\bf E_0},
\]
by Proposition \ref{prop:second-dual}. $\blacksquare$

\begin{Corollary} The operator ${\bf E_0}$ is injective and the operator ${\bf D_0}$ is surjective. Furthermore, the spaces ${\rm Ker}\, {\bf D_0}$ and ${\rm Im}\, {\bf E_0}$ are closed in $L^2_{\rm sym}(\Omega)$, and they are orthogonal complements of each other. 
\label{cor:E0-D0-properties}
\end{Corollary}
\noindent {\bf Proof.} From Theorem \ref{th:auxiliary_operators_are_mutually_adjoint} and Proposition \ref{prop:closed-operator-closed-kernel} it follows that ${\rm Ker}\, {\bf D_0}$ is closed, and by Proposition \ref{prop:prelim-kernel-of-the-adjoint}
\[
{\rm Ker}\, {\bf D_0} = \left({\rm Im}\, {\bf E_0}\right)^\perp.
\]
Next, recall the operator ${\rm E}$ defined by \eqref{eq:E-def1}--\eqref{eq:E-def2}, and notice that for some $c>0$ 
\[
\|u\|_{L^{2}}\leqslant  \|u\|_{W^{1,2}}\leqslant  c\|\nabla u\|_{L^2} \leqslant c\sqrt{2}\,\left\|{\rm E}(u)\right\|_{L^2}= c\sqrt{2}\,\left\|{\bf E_0}\, u\right\|_{L^2}\qquad \text{for any }u\in [W^{1,2}_0(\Omega)]^n,
\]
by, respectively, definition \eqref{eq:displacement_norms}, the Poincar\'{e}-Friedrichs inequality (see e. g. \cite[Th. 13.19, p.~430]{Leoni2017}, \cite[Th. 2.2, p.~74]{Necas1981}), Korn's first inequality (Proposition \ref{prop:Korn1}) and definition \eqref{eq:E0-def2}. We have fulfilled the condition of Proposition \ref{prop:Brezis-surjectivity} {\it ii)} with $A^*= {\bf E_0}$ defined as \eqref{eq:E0-def1}--\eqref{eq:E0-def2}. Thus, by Proposition \ref{prop:Brezis-surjectivity} and \eqref{eq:E0-D0-mutually-adjoint}, we have ${\rm Im} \,{\bf E_0}$ closed, ${\rm Ker}\, {\bf E_0}=\{0\}$ (i.e. ${\bf E_0}$ is injective) and ${\bf D_0}$ surjective. From Proposition \ref{prop:prelim-image-of-the-adjoint} and the closedness of ${\rm Im} \,{\bf E_0}$ we deduce the orthogonality
\[
{\rm Im}\, {\bf E_0}=\left({\rm Ker}\, {\bf D_0}\right)^\perp.
\]
$\blacksquare$

\subsection{The modified pair of adjoint operators accounting for the elasticity tensor  in the case of displacement boundary conditions}
\label{ssect:linear_problem1_operators_with_C}
Before we can solve the stress problem using the pair of adjoint operators ${\bf E_0, D_0}$ it is necessary to account for elasticity tensor $\bf C$ being non-constant with respect to $x$, see Section \ref{sssect:stresses}. Below we modify the pair of adjoint operators for this purpose.

First, consider space $L^{2}_{{\rm sym}, {\bf C}^{-1}}(\Omega)$, which is $L^2_{\rm sym}(\Omega)$, equipped with the following inner product and norm:
\begin{gather}
\la\tau_1, \tau_2\ra_{{\bf C}^{-1}} = \intOm \tau_1 : {\bf C}^{-1} \tau_2\, dx,\qquad \text{for any }\tau_1, \tau_2 \in L^2_{\rm sym}(\Omega), \label{eq:weighted-ip} \\
\|\tau\|_{{\bf C}^{-1}}=\left(\la\tau, \tau\ra_{{\bf C}^{-1}}\right)^{\frac{1}{2}}=\left(\intOm \tau : {\bf C}^{-1} \tau \, dx\right)^{\frac{1}{2}} \qquad \text{for any }\tau \in L^2_{\rm sym}(\Omega). \label{eq:weighted-norm}
\end{gather}
We will use the symbol $\perp_{{\bf C}^{-1}}$ to denote orthogonality in $L^{2}_{{\rm sym}, {\bf C}^{-1}}(\Omega)$, i.e. with respect to inner product \eqref{eq:weighted-ip}.
Since \eqref{eq:bounded_and_coercive_C_inv}--\eqref{eq:symmetric_C_inv} hold, form \eqref{eq:weighted-ip} is indeed an inner product, and norm \eqref{eq:weighted-norm} is equivalent to the usual norm $\|\cdot\|_{L^2}$. 
We define the operators
\[
{\bf E_{0, C}}:D({\bf E_{0,C}}) =[W_0^{1,2}(\Omega)]^n \subset  [L^2(\Omega)]^n \to L^2_{{\rm sym}, {\bf C}^{-1}}(\Omega),
\]
\[
{\bf E_{0, C}}: u \mapsto {\bf C}\, {\bf E_0}\,u
\]
\[
{\bf D_{0, C}}:D({\bf D_{0,C}}) = H_{\rm sym}({\rm div}, \Omega) \subset  L^2_{{\rm sym}, {\bf C}^{-1}}(\Omega)\to [L^2(\Omega)]^n,
\]
\begin{equation}
{\bf D_{0, C}}: \sigma \mapsto {\bf D_{0}}\, \sigma,
\label{eq:D0C-def2}
\end{equation}
for which Theorem \ref{th:auxiliary_operators_are_mutually_adjoint} and Corollary \ref{cor:E0-D0-properties} yield the following fact.
\begin{Corollary} Both operators ${\bf E_{0, C}}$, ${\bf D_{0,C}}$ are unbounded, closed, densely defined, and mutually adjoint. Moreover, the operator ${\bf E_{0, C}}$ is injective and the operator ${\bf D_{0, C}}$ is surjective. In turn, the spaces ${\rm Ker}\, {\bf D_{0,C}}$ and ${\rm Im}\, {\bf E_{0,C}}$ are closed in $L^2_{{\rm sym}, {\bf C}^{-1}}(\Omega)$, and they are orthogonal complements of each other:
\[
\left({\rm Ker}\, {\bf D_{0, C}}\right)^{\perp_{{\bf C}^{-1}}} = {\rm Im}\, {\bf E_{0, C}}.
\]
\end{Corollary}

\subsection{Solution to the static stress problem of linear elasticity with the displacement boundary conditions}
\label{ssect:linear_problem1_solution}
By applying Proposition \ref{prop:inverse_on_U} and Theorem \ref{th:abstract_linear_framework_solution} with 
\[
\mathcal{H} = L^{2}_{{\rm sym}, {\bf C}^{-1}}(\Omega),\qquad \mathcal{U} = {\rm Im}\, {\bf E_{0,C}}, \qquad \mathcal{V} = {\rm Ker}\, {\bf D_{0, C}},\qquad \mathcal{X}= [L^2(\Omega)]^n, \qquad A =  {\bf D_{0, C}}
\]
we obtain the following statement.
\begin{Theorem}
For any $g\in L^{2}_{{\rm sym}, {\bf C}^{-1}}(\Omega)$ and $f\in [L^2(\Omega)]^n$ the system with the unknown $\sigma\in  L^{2}_{{\rm sym}, {\bf C}^{-1}}(\Omega)$
\begin{numcases}{}
\sigma\in {\rm Im}\, {\bf E_{0, C}}+g, \label{eq:linear_problem1_eq1}\\
{\bf D_{0, C}}\,\sigma=f \label{eq:linear_problem1_eq2}
\end{numcases}
is equivalent to finding
\begin{equation}
\sigma\in ({\rm Im}\, {\bf E_{0, C}}+g) \cap( {\rm Ker}\, {\bf D_{0, C}}+Qf), 
\label{eq:intersection1}
\end{equation}
in which 
\[
{\rm Im}\, {\bf E_{0,C}}+g = {\rm Im}\, {\bf E_{0,C}}+P_{\,\mathcal{V}}\,g.
\]
Here $Q$ is the inverse of the restriction of ${\bf D_{0, C}}$ to ${{\rm Im}\, {\bf E_{0,C}}}$ and $P_{\,\mathcal{V}}$ is the orthogonal (in the sense of \eqref{eq:weighted-ip}) projection onto ${\rm Ker}\, {\bf D_{0, C}}$:
\[
Q:[L^2(\Omega)]^n \to {\rm Im}\, {\bf E_{0,C}}, \qquad P_{\,\mathcal{V}}: L^{2}_{{\rm sym}, {\bf C}^{-1}}(\Omega) \to {\rm Ker}\, {\bf D_{0, C}},
\]
and both $Q$ and $P_{\,\mathcal{V}}$ are continuous linear operators (their ranges have norm \eqref{eq:weighted-norm}).
Inclusion \eqref{eq:intersection1} has exactly one solution 
\begin{equation}
\sigma=P_{\,\mathcal{V}}g+Qf.
\label{eq:linear_solution1}
\end{equation} 
\label{th:linear_problem1_solution}
\end{Theorem}
Now, we confirm that the system \eqref{eq:linear_problem1_eq1}--\eqref{eq:linear_problem1_eq2} is indeed another way of writing a linear elasticity problem.
\begin{Proposition}Consider the {\it static stress problem of linear elasticity with the displacement boundary conditions} defined by equations  \eqref{eq:le-1'}, \eqref{eq:le-2}--\eqref{eq:le-3}, \eqref{eq:le-4'}, see Definitions \ref{def:static-problem-of-linear-elasticity-displacement} and \ref{def:reduced-problem-of-linear-elasticity}.  Assigning
\begin{equation}
f=F,\quad g={\bf C}\, {\rm E}(u_0),
\label{eq:loads_to_plug1}
\end{equation}
where  the operator ${\rm E}$ defined as \eqref{eq:E-def1}--\eqref{eq:E-def2}, we have
equation \eqref{eq:linear_problem1_eq1} equivalent to \eqref{eq:le-1'}, \eqref{eq:le-2}--\eqref{eq:le-3}, and equation \eqref{eq:linear_problem1_eq2} equivalent to \eqref{eq:le-4'}. 

\label{prop:linear_problem1_with_loads}
\end{Proposition}
\noindent {\bf Proof.} Let $\sigma$ satisfy \eqref{eq:linear_problem1_eq1}. Observe that 
\[
g={\bf C}\, {\rm E}(u_0)\in L^{2}_{\rm sym}(\Omega)\cong L^{2}_{{\rm sym}, {\bf C}^{-1}}(\Omega),
\]
and by plugging it to \eqref{eq:linear_problem1_eq1}, we obtain 
\begin{equation}
\sigma= {\bf C}\, {\rm E}(w) + {\bf C}\, {\rm E}(u_0)
\label{eq:in-proof-all-together1}
\end{equation}
with some $w\in [W^{1,2}_{0}(\Omega)]^n=D({\bf E_0})$. Take $u:=w+u_0$ and $\varepsilon:={\bf C}^{-1}\sigma$. Clearly, for the triple $(u, \varepsilon, \sigma)$ equations \eqref{eq:le-1'} and \eqref{eq:le-3} hold. Moreover, from \eqref{eq:in-proof-all-together1} we have that
\[
\varepsilon = {\rm E}(w+u_0),
\]
i.e. equation \eqref{eq:le-2} holds as well. 

Let $\sigma$ satisfy \eqref{eq:linear_problem1_eq2}. With $f=F$ and  \eqref{eq:D0C-def2},\, \eqref{eq:D0-def2} plugged in, it becomes
\[
-{\rm div}\, \sigma=F,
\]
which, by Definition \ref{def:H_div_sym-def} and Proposition \ref{prop:E-grad-symmetry-equality}, means that  
\[
\intOm \sigma:{\rm E}(u)\, dx = \intOm F\cdot u \, dx \qquad \text{for any }u\in [W^{1,2}_0(\Omega)]^n,
\]
which is exactly \eqref{eq:le-4'}. Vice versa, \eqref{eq:le-4'} implies \eqref{eq:linear_problem1_eq2} by the same reasoning.

Let \eqref{eq:le-1'}, \eqref{eq:le-2}--\eqref{eq:le-3} hold. By applying ${\bf C}\, {\rm E}$ to \eqref{eq:le-1'}, we get that
\[
{\bf C}\, {\rm E} (u) \in {\bf C}\, {\rm E}\left([W^{1,2}_{0}(\Omega)]^n\right)+ {\bf C}\, {\rm E}(u_0).
\]
Plug  \eqref{eq:le-2}--\eqref{eq:le-3} to the left-hand side and \eqref{eq:E0-def1}--\eqref{eq:E0-def2} to the right-hand side to obtain \eqref{eq:linear_problem1_eq1}. $\blacksquare$
\begin{Remark} Observe that by the continuity of $P_{\,\mathcal{V}}, Q, {\bf C}$ (in the sense of \eqref{eq:C-def3}) and ${\rm E}$ (in the sense of \eqref{eq:E-def1}) linear mapping \eqref{eq:linear_solution1}--\eqref{eq:loads_to_plug1} acting as
\begin{gather*}
[W^{1,2}(\Omega)]^n \times [L^{2}(\Omega)]^n \to L^{2}_{\rm sym}(\Omega)\cong L^2_{{\rm sym}, {\bf C}^{-1}}(\Omega),
\\
(u_0, F)\mapsto \sigma
\end{gather*}
is continuous. Moreover, Theorem \ref{th:linear_problem1_solution} and Proposition \ref{prop:linear_problem1_with_loads} mean that the schematic representation of Fig. \ref{fig:spaces} is suitable for the static stress problem of linear elasticity with the displacement boundary conditions.
\end{Remark}
\begin{Remark}
\label{remark:explicit-representations-displacement-bc}
For the problem with the displacement boundary conditions we can express the space $\mathcal{V}$ and the operators $P_{\,\mathcal{V}}$ and $Q$ in a more explicit way. 

\noindent\underline{Space $\mathcal{V}$ written explicitly.} From \eqref{eq:D0-def1}--\eqref{eq:D0-def2} we get that
\[
\mathcal{V} = {\rm Ker}\, {\bf D_{0, C}}  =\left\{\sigma\in H_{\rm sym}({\rm div}, \Omega): {\rm div}\, \sigma = 0\right\}.
\]

\noindent\underline{Operator $P_{\,\mathcal{V}}$ written explicitly.}
Let 
\[
P_{\,\mathcal{V}}\, g=\sigma
\]
for some $g, \sigma\in L^{2}_{{\rm sym}, {\bf C}^{-1}}(\Omega)$, which means  \eqref{eq:projection_abstract_def} with inner product \eqref{eq:weighted-ip}, i.e.
\[
\begin{cases} {\bf D_{0, C}},\sigma=0, \\[1mm]
\intOm\tau: {\bf C}^{-1}\sigma\, dx =\intOm \tau:{\bf C}^{-1} g\, dx & \text{for any }\tau \in L^{2}_{{\rm sym}, {\bf C}^{-1}}(\Omega)\text{ such that } {\bf D_{0, C}}\, \tau=0,
\end{cases}
\]
or
\[
\begin{cases} 
\sigma\in H_{\rm sym}({\rm div}, \Omega),\\[1mm]
{\rm div}\,\sigma=0, \\[1mm]
\intOm\tau: {\bf C}^{-1}\sigma\, dx =\intOm \tau:{\bf C}^{-1} g\, dx & \text{for any }\tau \in H_{\rm sym}({\rm div}, \Omega) \text{ with } {\rm div}\, \tau=0.
\end{cases}
\]

\noindent\underline{Operator $Q$ written explicitly.} Let
\[
Q\,f=\sigma
\]
for some $f\in [L^{2}(\Omega)]^n, \sigma\in  L^{2}_{{\rm sym}, {\bf C}^{-1}}(\Omega)$. By the construction of $Q$ it means that
\[
\begin{cases}
\sigma\in {\rm Im}\, {\bf E_{0,C}} = \left({\rm Ker}\, {\bf D_{0, C}}\right)^{\perp_{{\bf C}^{-1}}},\\
{\bf D_{0, C}}\,\sigma=f,
\end{cases}
\]
i.e.
\[
\begin{cases}
\intOm \tau: {\bf C}^{-1}\sigma=0& \text{for any }\tau \in H_{\rm sym}({\rm div}, \Omega) \text{ with } {\rm div}\, \tau=0,\\[1mm]
\sigma\in H_{\rm sym}({\rm div}, \Omega),\\
-{\rm div}\, \sigma=f.
\end{cases}
\]
\end{Remark}

\section{Analysis of the static stress problem of linear elasticity with the mixed boundary conditions}
\label{sect:problem2}
\subsection{Preliminaries: Lions-Magenes spaces and traces on the parts of the boundary}
\label{ssect:L-M_spaces}
The mixed boundary conditions lead to an issue of defining an appropriate space of functions for the traces of functions from $W^{1,2}(\Omega)$ and $H({\rm div}, \Omega)$ {\it on a part of the boundary}, in our case $\Gamma_{\rm D}$ or $\Gamma_{\rm T}$. Specifically, recall, that in \eqref{eq:le-4} we have $\delta u\in [W_{\rm D}]^n$, i.e. ${\rm tr (\delta u)}=0$ a.e. on $\Gamma_{\rm D}$, so that the last integral in \eqref{eq:le-4} is actually a linear functional over the {\it restrictions to $\Gamma_{\rm T}$ of traces of functions from $[W_{\rm D}]^n$}:
\[
\widetilde{T}: v \mapsto \int\limits_{\Gamma_{\rm T}} T\cdot v\, dS,
\]
\[
\widetilde{T}: D\left(\widetilde{T}\right)=\left\{{\rm tr}(\delta u)|_{\Gamma_{\rm T}}:\delta u\in [W^{1,2}(\Omega)]^n, {\rm tr}(\delta u)|_{\Gamma_{\rm D}}=0\right\} \to \mathbb{R}.
\]
For each $v\in D\left(\widetilde{T}\right)$ and the corresponding $\delta u$ we have, due to $\Gamma_{\rm T}\subset \Gamma$, that
\[
\|v\|_{[H^{1/2}(\Gamma_{\rm T})]^n}\leqslant \|{\rm tr} (\delta u)\|_{[H^{1/2}(\Gamma)]^n}
\]
with norm \eqref{eq:H1/2-norm-def}, therefore $v\in [H^{1/2}(\Gamma_{\rm T})]^n$. We cannot expect, however, that such functions $v$ will cover the entire space $[H^{1/2}(\Gamma_{\rm T})]^n$. 

Similarly, observe that one can take an arbitrary function $v\in H^{1/2}(\Gamma_{\rm T})$ and complete it by zero values on $\Gamma_{\rm D}$ to obtain a function $\widetilde{v}:\Gamma\to \mathbb{R}$. However, we cannot expect that $\widetilde{v}\in H^{1/2}(\Gamma)$ in general, as it would mean certain regularity throughout $\Gamma$, which is destroyed by the rough operation of ``gluing'' between $v$ and the constant-zero function (see also \cite[Remark 2.1.2, p.~49]{Boffi2013}).

The appropriate functional spaces for such traces are {\it Lions-Magenes spaces}, which can be defined as the via the following Proposition \ref{prop:zero-trace-on-GD-functions} and Definition \ref{def:L-M-spaces}.
\begin{Proposition}
\label{prop:zero-trace-on-GD-functions}
Recall that $\Gamma_{\rm T}\subset \Gamma$ is a nonempty subset, open in $\Gamma$. The following subspace of $H^{1/2}(\Gamma)$
\begin{equation}
H^{1/2}_{\rm D}(\Gamma) = \{u\in H^{1/2}(\Gamma): u|_{\Gamma_{\rm D}}=0 \text{ a.e.}\}
\label{eq:zero-trace-on-GD-functions-def}
\end{equation}
is closed.
\end{Proposition}
\noindent {\bf Proof.} Consider $\{u_k\}_{k\in\mathbb{N}}\subset H^{1/2}_{\rm D}(\Gamma)$ and $u^*\in H^{1/2}(\Gamma)$ such that $u_k \to u^*$ in $H^{1/2}(\Gamma)$. Then
\[
\left|\,\int\limits_{\Gamma_{\rm D}} \left(u_k-u^*\right)^2\,dx\right|\leqslant \|u_k-u^*\|^2_{L^2(\Gamma)}\leqslant \|u_k-u^*\|^2_{H^{1/2}(\Gamma)} \to 0,
\]
see \eqref{eq:H1/2-norm-def} for the last inequality. Therefore, $u^*$ is zero a.e. on $\Gamma_{\rm D}$ (see e.g. \cite[Th. 4.9, p.~94]{Brezis2011}), from which the closedness of $H^{1/2}_{\rm D}(\Gamma)$ follows. $\blacksquare$

\begin{Definition}{\bf} \cite[pp.~36--37]{Gatica2014}, \cite[p.~121]{Sayas2019}, \cite[pp.~35--36]{Chouly2023}
The {\it Lions-Magenes space} is
\begin{equation}
H^{1/2}_{00} (\Gamma_{\rm T})=\left\{u|_{\Gamma_{\rm T}}: u\in H^{1/2}_{\rm D}(\Gamma)\right\}
\label{eq:L-M-space-1-def}
\end{equation}
with the norm
\begin{equation}
\|u|_{\Gamma_{\rm T}}\|_{H^{1/2}_{00} (\Gamma_{\rm T})}= \|u\|_{H^{1/2} (\Gamma)} \qquad \text{ for }u\in H_{\rm D}^{1/2}(\Gamma).
\label{eq:L-M-space-norm-def}
\end{equation}
We denote the  {\it dual Lions-Magenes space} 
\[
H^{-1/2}_{00} (\Gamma_{\rm T})= H^{1/2}_{00} (\Gamma_{\rm T})^*
\]
and the restriction of the trace operator
\begin{gather}
{\rm tr\,}_{\Gamma_{\rm T}}: W_{\rm D} \to H^{1/2}_{00} (\Gamma_{\rm T}),
\label{eq:tr-restriction-def1}\\[2mm]
{\rm tr\,}_{\Gamma_{\rm T}}: u \mapsto {\rm tr}(u)|_{\Gamma_{\rm T}}.
\label{eq:tr-restriction-def2}
\end{gather}
\label{def:L-M-spaces}
\end{Definition}
\begin{Remark} Observe that ${\rm tr\, }_{\Gamma_{\rm T}}$ is surjective, due to \eqref{eq:L-M-space-1-def} and the surjectivity of the trace operator to $H^{1/2}(\Gamma)$, see Theorem \ref{th:trace_th_1}. Also recall that $W_{\rm D}$ is equipped with the norm of $W^{1,2}(\Omega)$ (see \eqref{eq:W-def}), so that the boundness of the trace operator (see Theorem \ref{th:trace_th_1} again) implies the boundness of ${\rm tr\,}_{\Gamma_{\rm T}}$ as defined in \eqref{eq:tr-restriction-def1}--\eqref{eq:tr-restriction-def2} since for $u\in W_{\rm D}$
\begin{equation}
\|{\rm tr\,}_{\Gamma_{\rm T}}(u)\|_{H^{1/2}_{00} (\Gamma_{\rm T})}=\|{\rm tr}(u)\|_{H^{1/2}(\Gamma)}\leqslant \|{\rm tr}\|_{\rm op}\, \|u\|_{W^{1,2}}.
\label{eq:tr-restriction-boundness}
\end{equation}
\label{remark:restriction_of_traces_surjectivity}
\end{Remark}
\begin{Proposition}
The space  $H^{1/2}_{00} (\Gamma_{\rm T})$ is Hilbert and, therefore, reflexive.
\label{prop:L-M_spaces_are_reflexive}
\end{Proposition}
\noindent {\bf Proof.} Observe that the space $H_{\rm D}^{1/2}(\Gamma)$ is a Hilbert space with inner product \eqref{eq:H1/2-ip-def} due to Proposition \ref{prop:zero-trace-on-GD-functions}.

Consider the linear operators of {\it restriction to $\Gamma_{\rm T}$} and {\it extension by zero}:
\[
{\rm r}: H^{1/2}_{\rm D}(\Gamma) \to H^{1/2}_{00} (\Gamma_{\rm T}), \qquad {\rm z}:H^{1/2}_{00} (\Gamma_{\rm T}) \to H^{1/2}_{\rm D}(\Gamma),
\]
\[
{\rm r}:u \mapsto u|_{\Gamma_{\rm T}}, \qquad {\rm z}: u \mapsto \left(\begin{cases}u(x),& x\in \Gamma_{\rm T}\\ 0, & x\in\Gamma_{\rm D}\end{cases}\right).
\]
They are mutually inverse (hence they are bijective) and norm-preserving (hence they are continuous). Therefore, the following form, defined on $H^{1/2}_{00} (\Gamma_{\rm T})$, is symmetric, linear and positive definite:
\begin{multline}
\la u,v\ra_{H^{1/2}_{00} (\Gamma_{\rm T})} := \la {\rm z}(u), {\rm z}(v)\ra_{H^{1/2}(\Gamma)}=\la u, v\ra_{H^{1/2}(\Gamma_{\rm T})}+2\int\limits_{\Gamma_{\rm T}} u(x)\,v(x) \int\limits_{\Gamma_{\rm D}}\frac{1}{\|x-y\|^n}dS_y dS_x \\ \text{for any}\, u,v\in H^{1/2}_{00}(\Gamma_{\rm T}).
\label{eq:L-M-space-ip-def}
\end{multline}
In \eqref{eq:L-M-space-ip-def} the last equality can be shown directly by the substitution of \eqref{eq:H1/2-ip-def}.
The corresponding norm to \eqref{eq:L-M-space-ip-def} is \eqref{eq:L-M-space-1-def}, so that the completeness of  $H^{1/2}_{\rm D}(\Gamma)$ and the properties of operators ${\rm r}$ and ${\rm z}$ imply that $H^{1/2}_{00} (\Gamma_{\rm T})$ is a Hilbert space with the inner product \eqref{eq:L-M-space-ip-def}. 
 $\blacksquare$

Some authors use notation $\widetilde{H}^{1/2}(\Gamma_{\rm T})$ for space $H^{1/2}_{00} (\Gamma_{\rm T})$ itself \cite[p.~121]{Sayas2019} or for $H_{\rm D}^{1/2}(\Gamma)$ obtained as $\overline{\mathcal{D}(\Gamma_{\rm T})}^{H^{1/2}(\Gamma)}$ \cite[p.~99]{McLean2000}. We also refer to \cite[p.~84]{Kikuchi1988}, \cite[p.~379]{Baiocchi1984} \cite[Section 33, pp.~159-161]{Tartar2007}, \cite[Th. 6.97, p.~245]{Leoni2023} and the original source \cite[Th. 11.7, p.~66, Remark 12.1, p.~71]{Lions1972} for different definitions of Lions-Magenes spaces and for an intrinsic definition of a norm in those spaces. Such intrinsic definitions, however, may require additional regularity of $\Gamma_{\rm T}$ to guarantee that the functions are extendable to $\Gamma$ and $\Omega$, see Remark \ref{remark:Lipschitz-dissection}.

Similarly to \eqref{eq:vector-valued-trace-def}, for $\mathbb{R}^n$-valued functions $u=(u_i)_{i\in \overline{1,n}} \in[W_{\rm D}]^n$ we denote
\[
{\rm tr \,}_{\Gamma_{\rm T}}(u):= ({\rm tr\,}_{\Gamma_{\rm T}}(u_i))_{i\in \overline{1,n}}\in \left[H^{1/2}_{00}(\Gamma_{\rm T})\right]^n.
\]

\begin{Theorem}{\bf (normal trace on a part of the boundary)} \cite[Proposition 2.19, p.~51]{Chouly2023}
In the setting of Definition \ref{def:L-M-spaces} there exists a continuous linear operator (called the normal trace operator)
\begin{equation}
{\rm tr}_{\nu, \Gamma_{\rm T}}: H({\rm div}, \Omega) \to \left[H^{-1/2}_{00} (\Gamma_{\rm T})\right]^n,
\label{eq:normal_trace_of_a_matrix_field}
\end{equation}
where $H({\rm div}, \Omega)\subset [L^2(\Omega)]^{n\times n}$ is as in Definition \ref{def:H_div_sym-def} for which
\[
\left \la {\rm tr}_{\nu, \Gamma_{\rm T}}(\tau), {\rm tr\,}_{\Gamma_{\rm T}}(u) \right\ra = \int\limits_{\Gamma_{\rm T}}u\cdot (\tau \nu)\, dS \qquad \text{ for any }\tau\in [\mathcal{D}(\overline{\Omega})]^{n\times n}, \, u\in [\mathcal{D}_{\rm D}]^n,
\]
where $\mathcal{D}_{\rm D}$ is as in \eqref{eq:D-pb-def} and  $\nu(x) \in \mathbb{R}^n$ is a unit outward normal which exists at a.a. $x\in \Gamma$. Such operator is unique.
\label{th:normal_trace_of_a_matrix_field}
\end{Theorem}
\begin{Theorem} {\bf (Green's formula for matrix-valued divergence on a part of the boundary)}
The following equalities hold:
\begin{equation}
\left \la {\rm tr}_{\nu, \Gamma_{\rm T}}(\tau), {\rm tr\,}_{\Gamma_{\rm T}}(u) \right\ra = \intOm \tau: \nabla u\, dx+  \intOm ({\rm div}\, \tau) \cdot u \, dx \qquad \text{for any } \tau\in H({\rm div}, \Omega), \, u\in [W_{\rm D}]^n
\label{eq:Green's_formula_for_matrices1}
\end{equation}
and 
\begin{equation}
\left \la {\rm tr}_{\nu, \Gamma_{\rm T}}(\tau), {\rm tr\,}_{\Gamma_{\rm T}}(u) \right\ra = \intOm \tau: {\rm E} (u)\, dx+  \intOm ({\rm div}\, \tau) \cdot u \, dx \qquad \text{for any } \tau\in H_{\rm sym}({\rm div}, \Omega), \, u\in [W_{\rm D}]^n,
\label{eq:Green's_formula_for_matrices2}
\end{equation}
where ${\rm E}$ is as in \eqref{eq:E-def1}--\eqref{eq:E-def2}.
\label{th:Green's_formula_for_matrices}
\end{Theorem}
\noindent {\bf Proof.} For the statement and the poof of Green's formula on the part of the boundary with scalar-valued $u$ we refer to \cite[Proposition 2.20, p.~52]{Chouly2023}. From it \eqref{eq:Green's_formula_for_matrices1} follows directly, and \eqref{eq:Green's_formula_for_matrices2} follows then from Proposition \ref{prop:E-grad-symmetry-equality}. $\blacksquare$

\noindent Regarding Theorems \ref{th:normal_trace_of_a_matrix_field}-\ref{th:Green's_formula_for_matrices} see also \cite[Section 3.2.3, pp.~72--73]{Chouly2023}, \cite[pp.~121--122]{Sayas2019}, \cite[p.~47]{Gatica2014} and \cite[Th. 5.6, p.~88, Th. 5.11, p.~93]{Kikuchi1988} (the latter has further decomposition of the trace spaces into normal and tangential components). Green's formula for scalar-valued functions can be found in \cite[(1.2.4), p.~14]{Ciarlet1978}, \cite[(18.1), p.~591]{Leoni2017}.

\subsection{The pair of adjoint operators accounting for traces}
\label{ssect:main_operators}
Now we consider the following operators 
\begin{gather}
{\bf D_1}:D\left({\bf D_1}\right) = H_{\rm sym}({\rm div}, \Omega) \subset  L^2_{\rm sym}(\Omega)\to [L^2(\Omega)]^n \times [H^{-1/2}_{00}(\Gamma_{\rm T})]^n, \label{eq:D1-def1}
\\
{\bf D_1}: \sigma\mapsto \begin{pmatrix} - {\rm div}\, \sigma \\  {\rm tr}_{\nu, \Gamma_{\rm T}} (\sigma)
\end{pmatrix},
\label{eq:D1-def2}
\end{gather}
\begin{gather}
{\bf E_1}: D\left({\bf E_1}\right) = \left\{(u, {\rm tr\,}_{\Gamma_{\rm T}}(u)): u\in [W_{\rm D}]^n\right\}\subset [L^2(\Omega)]^n\times [H^{1/2}_{00}(\Gamma_{\rm T})]^n \to L^2_{\rm sym}(\Omega), \label{eq:E1-def1}
\\
{\bf E_1}:\begin{pmatrix} u\\ {\rm tr\,}_{\Gamma_{\rm T}}(u)\end{pmatrix} \mapsto \frac{1}{2}\left(\nabla u+ (\nabla u)^{\top}\right).
\label{eq:E1-def2}
\end{gather}

\begin{Theorem}
\label{th:main-pair}
Both ${\bf D_1}$ and ${\bf E_1}$ are closed densely defined unbounded operators in the sense of Proposition \ref{prop:adjoint-def} and 
\begin{equation}
{\bf D_1}^* = {\bf E_1} \quad \text{ and }\quad {\bf E_1}^* = {\bf D_1}.
\label{eq:E1-D1-mutually-adjoint}
\end{equation}
In particular, the density of the domain of ${\bf E_1}$ means that
\begin{equation}
\overline{\left\{\vphantom{\widetilde{V}} (u, {\rm tr\,}_{\Gamma_{\rm T}}(u)): u\in [W_{\rm D}]^n\right\}}^{\,[L^2(\Omega)]^n\times [H^{1/2}_{00}(\Gamma_{\rm T})]^n}= [L^2(\Omega)]^n\times [H^{1/2}_{00}(\Gamma_{\rm T})]^n.
\label{eq:independent-density-statement}
\end{equation}
\end{Theorem}
\noindent {\bf Proof.} \newline
\underline{Claim: ${\bf D_1}$ is densely defined.}
From Theorem \ref{th:auxiliary_operators_are_mutually_adjoint} (specifically, \eqref{eq:density_H_sym_in_L2}) we already know that ${\bf D_1}$ is densely defined.
\newline \underline{Claim: ${\bf D_1}$ is closed.} Let
\[
\{\sigma_k\}_{k\in \mathbb{N}}\subset H_{\rm sym}({\rm div}, \Omega), \quad \sigma^*\in L^2_{\rm sym}(\Omega), \quad (f, \eta) \in [L^2(\Omega)]^n \times [H^{-1/2}_{00}(\Gamma_{\rm T})]^n
\] 
be such that
\begin{gather}
\|\sigma_k-\sigma^*\|_{L^2} \to 0, \label{eq:adjoint-proof-L1}\\
\|-{\rm div}\, \sigma_k - f\|_{L^2} \to 0, \label{eq:adjoint-proof-L2}\\
\left\| {\rm tr}_{\nu, \Gamma_{\rm T}} (\sigma_k) - \eta \right\|_{[H^{-1/2}_{00}(\Gamma_{\rm T})]^n} \to 0. \label{eq:adjoint-proof-L3}
\end{gather}
To establish the closedness we need to show that:
\begin{gather}
\sigma^*\in H_{\rm sym}({\rm div}, \Omega), \label{eq:adjoint-proof-NTS1}\\
f= -{\rm div}\, \sigma^*,\label{eq:adjoint-proof-NTS2}\\
\eta = {\rm tr}_{\nu, \Gamma_{\rm T}} (\sigma^*). \label{eq:adjoint-proof-NTS3}
\end{gather} 
Observe that for any $u\in [\mathcal{D}(\Omega)]^n$
\begin{align*}
\left|\intOm \sigma^* : \nabla u\, dx -\intOm f\cdot u\, dx\right| &= \left|\intOm (\sigma^*-\sigma_k) : \nabla u\, dx +\intOm \sigma_k:\nabla u\, dx- \intOm f\cdot u\, dx\right|=\\
&=\left|\intOm (\sigma^*-\sigma_k) : \nabla u\, dx + \intOm (-f-{\rm div}\, \sigma_k)\cdot u\, dx\right|\leqslant \\[1mm]
&\leqslant \|\sigma^*-\sigma_k\|_{L^2}\, \|\nabla u\|_{L^2} +\|-f-{\rm div}\, \sigma_k\|_{L^2}\,\|u\|_{L^2} \to 0,
\end{align*}
where the zero limit is due to \eqref{eq:adjoint-proof-L1} and \eqref{eq:adjoint-proof-L2}. Thus, by Definition \ref{def:H_div_sym-def}  we have $\sigma^*\in  H_{\rm sym}({\rm div}, \Omega)$ with ${\rm div}\, \sigma^* =-f$. We have shown \eqref{eq:adjoint-proof-NTS1} and \eqref{eq:adjoint-proof-NTS2}. Moreover, \eqref{eq:adjoint-proof-L1}, \eqref{eq:adjoint-proof-L2} and \eqref{eq:adjoint-proof-NTS2} together  mean that
\begin{equation}
\|\sigma_k-\sigma^*\|_{H({\rm div}, \Omega)} \to 0,
\label{eq:H-div-convergence-in-the-proof}
\end{equation}
where the norm was defined as \eqref{eq:H-div-norm}. Together limits \eqref{eq:H-div-convergence-in-the-proof} and \eqref{eq:adjoint-proof-L3} imply \eqref{eq:adjoint-proof-NTS3} due to the continuity of normal trace operator \eqref{eq:normal_trace_of_a_matrix_field} and the uniqueness of the limit. We have proven that ${\bf D_1}$ is closed.
\newline \underline{Claim: $D\left({\bf D_1}^*\right)\subset D\left({\bf E_1}\right)$.}
As defined by \eqref{eq:D(A)-def},
\begin{multline}
D\left({\bf D_1}^*\right) = \left\{\left(u, \widetilde{v}\right)\in [L^{2}(\Omega)^n]\times  \left(\left[H^{-1/2}_{00} (\Gamma_{\rm T})\right]^n\right)^*: \text{there exisits }c_1\geqslant 0: \vphantom{ \left|\intOm \right|}\right.\\ \left. \text{for all }\sigma \in H_{\rm sym}({\rm div}, \Omega)\quad \left|\intOm u\cdot (- {\rm div}\, \sigma)\, dx + \left\la\widetilde{v}, {\rm tr}_{\nu, \Gamma_{\rm T}}(\sigma)\right\ra \right|\leqslant c_1\|\sigma\|_{L^2}    \right\}=\\
=\left\{\left(u, {\rm tr\,}_{\Gamma_{\rm T}}(v)\right): u\in [L^{2}(\Omega)^n], v\in [W_{\rm D}]^n: \text{there exisits }c_1\geqslant 0: \vphantom{ \left|\intOm \right|}\right.\\ \left. \text{for all }\sigma \in H_{\rm sym}({\rm div}, \Omega)\quad \left|\intOm u\cdot (- {\rm div}\, \sigma)\, dx + \left\la{\rm tr}_{\nu, \Gamma_{\rm T}}(\sigma),{\rm tr\,}_{\Gamma_{\rm T}}(v) \right\ra \right|\leqslant c_1\|\sigma\|_{L^2}\right\},
\label{eq:D(-Div)-def}
\end{multline}
where the second equality follows from the reflexivity of $H^{1/2}_{00} (\Gamma_{\rm T})$ (see Proposition \ref{prop:L-M_spaces_are_reflexive}) and the surjectivity of ${\rm tr\,}_{\Gamma_{\rm T}}$ (see Remark \ref{remark:restriction_of_traces_surjectivity}). By Green's formula \eqref{eq:Green's_formula_for_matrices2} we have
\begin{align*}
\left|\intOm u\cdot (- {\rm div}\, \sigma)\, dx + \left\la{\rm tr}_{\nu, \Gamma_{\rm T}}(\sigma),{\rm tr\,}_{\Gamma_{\rm T}}(v) \right\ra \right|
&=\left|\intOm u\cdot (- {\rm div}\, \sigma)\, dx + \intOm \sigma: {\rm E} (v)\, dx+  \intOm ({\rm div}\, \sigma) \cdot v \, dx \right| =\\
&=\left|\intOm (u-v)\cdot (- {\rm div}\, \sigma)\, dx - \left(- \intOm \sigma: {\rm E} (v)\, dx\right) \right|\geqslant\\
&\geqslant \left|\intOm (u-v)\cdot (- {\rm div}\, \sigma)\, dx\right| - \left|\left(-\intOm \sigma: {\rm E} (v)\, dx \right)\right|.
\end{align*}
Thus for any $\sigma\in H_{\rm sym}({\rm div}, \Omega)$ and $u\in [L^{2}(\Omega)^n], \,v\in [W_{\rm D}]^n$ such that $\left(u, {\rm tr\,}_{\Gamma_{\rm T}}(v)\right)\in D\left({\bf D_1}^*\right)$ we have by \eqref{eq:D(-Div)-def} that
\[
\left|\intOm (u-v)\cdot (- {\rm div}\, \sigma)\, dx\right|\leqslant c_1\|\sigma\|_{L^2}+ \left|\left(-\intOm \sigma: {\rm E} (v)\, dx \right)\right|\leqslant c_1\|\sigma\|_{L^2}+ \|{\rm E}(v)\|_{L^2}\|\sigma\|_{L^2},
\]
i.e.
\[
\left| \intOm(u-v)\cdot ({\bf D_0}\, \sigma)\, dx \right|\leqslant c_2\|\sigma\|_{L^2},
\]
where operator ${\bf D_0}$ is discussed in Section \ref{ssect:aux_operators} and we set $c_2=c_1+\|{\rm E}(v)\|_{L^2}$. Observe that, since $c_2$ is independent of $\sigma$, 
\[
u-v \in D\left({\bf D_0}^*\right) = D({\bf E_0}) =  [W_0^{1,2}(\Omega)]^n
\]
by \eqref{eq:D(A)-def} and Theorem \ref{th:auxiliary_operators_are_mutually_adjoint}. Hence,
\[
u\in [W_0^{1,2}(\Omega)]^n + [W_{\rm D}]^n \subset [W^{1,2}(\Omega)]^n \qquad \text{ and } \qquad {\rm tr}(u) ={\rm tr}(v),
\]
the latter due to Proposition \ref{prop:W0-ker-tr}. Thus $u\in [W_{\rm D}]^n$ and we have
\[
D\left({\bf D_1}^*\right)\subset  \left\{(u, {\rm tr\,}_{\Gamma_{\rm T}}(u)): u\in [W_{\rm D}]^n\right\}.
\]
\noindent \underline{Finishing the proof.} Vice versa, consider $(u, {\rm tr\,}_{\Gamma_{\rm T}}(u))$ for an arbitrary $u\in [W_{\rm D}]^n$. Then, utilizing \eqref{eq:Green's_formula_for_matrices2} again, we have for any $\sigma\in L^2_{\rm sym}(\Omega)$ that
\begin{equation}
\begin{aligned}
\left\la \begin{pmatrix} u\\ {\rm tr\,}_{\Gamma_{\rm T}}(u)\end{pmatrix},{\bf D_1}\, \sigma\right\ra &= 
\intOm u\cdot (- {\rm div}\, \sigma)\, dx + \left\la{\rm tr}_{\nu, \Gamma_{\rm T}}(\sigma),{\rm tr\,}_{\Gamma_{\rm T}}(u) \right\ra =\\
&=\intOm {\rm E}(u):\sigma\, dx = \left\la {\bf E_1}\begin{pmatrix} u\\ {\rm tr\,}_{\Gamma_{\rm T}}(u)\end{pmatrix}, \sigma \right\ra
\end{aligned}
\label{eq:main_adjoint_relation_proof}
\end{equation}
and thus
\[
\left|\left\la \begin{pmatrix} u\\ {\rm tr\,}_{\Gamma_{\rm T}}(u)\end{pmatrix},{\bf D_1}\, \sigma\right\ra\right| = \left|\intOm {\rm E}(u):\sigma\, dx \right| \leqslant \|{\rm E}(u)\|_{L^2} \|\sigma\|_{L^2},
\]
i.e. 
\[
(u, {\rm tr\,}_{\Gamma_{\rm T}}(u))\in D\left({\bf D_1}^*\right) 
\]
according to \eqref{eq:D(-Div)-def}. We have proven that, indeed, 
\[
D\left({\bf D_1}^*\right)=\left\{(u, {\rm tr\,}_{\Gamma_{\rm T}}(u)): u\in [W_{\rm D}]^n\right\}=D\left({\bf E_1}\right).
\]
Then relation \eqref{eq:main_adjoint_relation_proof} means that
\[
{\bf D_1}^* = {\bf E_1} 
\]
and, therefore, ${\bf E_1}$ is a closed operator (see Proposition \ref{prop:adjoint-def}). Moreover, since both $L^2_{\rm sym}(\Omega)$ and $[L^2(\Omega)]^n \times  [H^{-1/2}_{00}(\Gamma_{\rm T})]^n$ are reflexive, $ {\bf E_1}$ is densely defined (i.e. \eqref{eq:independent-density-statement} holds) and 
\[
{{\bf E_1}}^* = {\bf D_1},
\]
by Proposition \ref{prop:second-dual}. $\blacksquare$
\begin{Corollary} The operator ${\bf E_1}$ is injective and the operator ${\bf D_1}$ is surjective.  Furthermore, the spaces ${\rm Ker}\, {\bf D_1}$ and ${\rm Im}\, {\bf E_1}$ are closed in $L^2_{\rm sym}(\Omega)$, and they are orthogonal complements of each other. 
\label{cor:E1-D1-properties}
\end{Corollary}
\noindent {\bf Proof.}
The proof is mostly identical to the proof of Corollary \ref{cor:E0-D0-properties}, but for operators ${\bf E_1, D_{1}}$ we rely on Theorem \ref{th:main-pair}. To fulfill the condition of Proposition \ref{prop:Brezis-surjectivity} {\it ii)} we use the following estimate. From the boundness of ${\rm tr}_{\Gamma_{\rm T}}$ (i.e. \eqref{eq:tr-restriction-boundness}) and the implication of {\it Korn's second inequality} \eqref{eq:fromKorn2} we have that for any $u\in [W_{\rm D}]^n$
\begin{multline*}
\left\|\begin{pmatrix} u\\ {\rm tr\,}_{\Gamma_{\rm T}}(u)\end{pmatrix}\right\|_{[L^2(\Omega)]^n\times [H^{1/2}_{00}(\Gamma_{\rm T})]^n} =\left( \|u\|_{L^2}^2+ \left\|{\rm tr\,}_{\Gamma_{\rm T}}(u)\right\|_{[H^{1/2}_{00}(\Gamma_{\rm T})]^n}^2\right)^{\frac{1}{2}}\leqslant\\[2mm]
\leqslant \left( \|u\|_{L^2}^2+ \left\|{\rm tr}\right\|_{\rm op}^2\,\|u\|_{W^{1,2}}^2\right)^{\frac{1}{2}}\leqslant c_1\|u\|_{W^{1,2}} \leqslant c_1c_2\,\|{\rm E}(u)\|_{L^2} = c_1c_2\, \left\|{\bf E_1}\begin{pmatrix} u\\ {\rm tr\,}_{\Gamma_{\rm T}}(u)\end{pmatrix}\right\|_{L^2}
\end{multline*}
with some constants $c_1,c_2\geqslant 0$.
$\blacksquare$

\begin{Remark}
In \cite[Section 5.1, Th. 5.1, pp.~3780--3781]{Amstutz2020} a different formalism was used to derive a decomposition of $L^2_{\rm sym}(\Omega)$, similar to ${\rm Ker}\, {\bf D_1} \oplus {\rm Im}\, {\bf E_1}$, under the assumption of a simply-connected domain. In there, the decomposition helps to establish a strain-based model of elastoplasticity, constructed via the incompatibility operator $\rm curl\, curl$. In the current text we focus on the formulation of elasticity in terms of stress and on the stress-to-forces mapping ${\bf D_1}$ in particular.
\end{Remark}

\subsection{The modified pair of adjoint operators accounting for the elasticity tensor in the case of mixed boundary conditions}
Again, we should modify the pair of adoint operators to account for ${\bf C}$ being non-constant with respect to $x$. Similarly to Section \ref{ssect:linear_problem1_operators_with_C}, we will use space $L^{2}_{{\rm sym}, {\bf C}^{-1}}(\Omega)$ taken as $L^{2}_{\rm sym}(\Omega)$ with inner product \eqref{eq:weighted-ip}, equivalent norm \eqref{eq:weighted-norm} and the orthogonality relation denoted by $\perp_{{\bf C}^{-1}}$. We define the operators
\begin{gather*}
{\bf E_{1, C}}: D\left({\bf E_1}\right) = \left\{(u, {\rm tr\,}_{\Gamma_{\rm T}}(u)): u\in [W_{\rm D}]^n\right\}\subset [L^2(\Omega)]^n\times [H^{1/2}_{00}(\Gamma_{\rm T})]^n \to L^2_{\rm sym, {\bf C}^{-1}}(\Omega), 
\\
{\bf E_{1,C}}:\begin{pmatrix} u\\ {\rm tr\,}_{\Gamma_{\rm T}}(u)\end{pmatrix} \mapsto  {\bf C}\, {\bf E_{1}} \begin{pmatrix} u\\ {\rm tr\,}_{\Gamma_{\rm T}}(u)\end{pmatrix},
\end{gather*}
\[
{\bf D_{1, C}}: D({\bf D_{1, C}}) = H_{\rm sym}({\rm div}, \Omega)\subset L^{2}_{{\rm sym},{\bf C}^{-1}}(\Omega)
\to [L^2(\Omega)]^n \times [H^{-1/2}_{00}(\Gamma_{\rm T})]^n,
\]
\begin{equation}
{\bf D_{1, C}}:u \mapsto {\bf D_{1}}\, u,
\label{eq:D1C-def2}
\end{equation}
for which Theorem \ref{th:main-pair} and Corollary \ref{cor:E1-D1-properties} yield the following fact.
\begin{Corollary} Both operators ${\bf E_{1, C}}$, ${\bf D_{1,C}}$ are unbounded, closed, densely defined, and mutually adjoint. Moreover, the operator ${\bf E_{1, C}}$ is injective and the operator ${\bf D_{1, C}}$ is surjective. In turn, the spaces ${\rm Ker}\, {\bf D_{1,C}}$ and ${\rm Im}\, {\bf E_{1,C}}$ are closed in $L^2_{{\rm sym}, {\bf C}^{-1}}(\Omega)$, and they are orthogonal complements of each other:
\[
\left({\rm Ker}\, {\bf D_{1, C}}\right)^{\perp_{{\bf C}^{-1}}} = {\rm Im}\, {\bf E_{1, C}}.
\]
\end{Corollary}
\subsection{Solution to the static stress problem of linear elasticity with the mixed boundary conditions}
\label{ssect:linear_problem2_solution}

Notice, that by Proposition \ref{prop:L-M_spaces_are_reflexive} the space $H^{-1/2}_{00}(\Gamma_{\rm T})$ is a Hilbert space, by the identification with $H^{1/2}_{00}(\Gamma_{\rm T})$ using the inner product \eqref{eq:L-M-space-ip-def}. Therefore, we can apply Proposition \ref{prop:inverse_on_U} and Theorem \ref{th:abstract_linear_framework_solution} with 
\[
\mathcal{H} = L^{2}_{{\rm sym}, {\bf C}^{-1}}(\Omega),\quad \mathcal{U} = {\rm Im}\, {\bf E_{1,C}}, \quad \mathcal{V} = {\rm Ker}\, {\bf D_{1, C}},
\]
\[
\mathcal{X}= [L^2(\Omega)]^n\times [H^{-1/2}_{00}(\Gamma_{\rm T})]^n, \quad A= {\bf D_{1, C}}
\]
to obtain the following statement.
\begin{Theorem}
For any $g\in L^{2}_{{\rm sym}, {\bf C}^{-1}}(\Omega),\, f\in [L^2(\Omega)]^n \times [H^{-1/2}_{00}(\Gamma_{\rm T})]^n$ the system with the unknown $\sigma\in  L^{2}_{{\rm sym}, {\bf C}^{-1}}(\Omega)$
\begin{numcases}{}
\sigma\in {\rm Im}\, {\bf E_{1,C}}+g, \label{eq:linear_problem2_eq1}\\
{\bf D_{1, C}}\,\sigma=f \label{eq:linear_problem2_eq2}
\end{numcases}
is equivalent to finding
\begin{equation}
\sigma\in ( {\rm Im}\, {\bf E_{1,C}}+g) \cap( {\rm Ker}\, {\bf D_{1, C}}+Qf), 
\label{eq:intersection2}
\end{equation}
in which
\[
{\rm Im}\, {\bf E_{1,C}}+g= {\rm Im}\, {\bf E_{1,C}}+P_{\,\mathcal{V}}\,g .
\]
Here $Q$ is the inverse of the restriction of ${\bf D_{1, C}}$ to ${\rm Im}\, {\bf E_{1,C}}$ and $P_{\,\mathcal{V}}$ is the orthogonal (in the sense of \eqref{eq:weighted-ip}) projection onto ${\rm Ker}\, {\bf D_{1, C}}$:
\[
Q:[L^2(\Omega)]^n\times [H^{-1/2}_{00}(\Gamma_{\rm T})]^n \to  {\rm Im}\, {\bf E_{1,C}}, \qquad P_{\,\mathcal{V}}: L^{2}_{{\rm sym}, {\bf C}^{-1}}(\Omega) \to {\rm Ker}\, {\bf D_{1, C}},
\]
and both $Q$ and $P_{\,\mathcal{V}}$ are continuous linear operators (their ranges have norm \eqref{eq:weighted-norm}).
Inclusion \eqref{eq:intersection2} has exactly one solution 
\begin{equation}
\sigma=P_{\,\mathcal{V}}g+Qf.
\label{eq:linear_solution2}
\end{equation} 
\label{th:linear_problem2_solution}
\end{Theorem}
Again, we verify the relation of the system \eqref{eq:linear_problem2_eq1}--\eqref{eq:linear_problem2_eq2} with a linear elasticity problem.

\begin{Proposition}Consider the {\it static stress problem of linear elasticity with the mixed boundary conditions} defined by equations  \eqref{eq:le-1}--\eqref{eq:le-4}, see Definitions \ref{def:static-problem-of-linear-elasticity-mixed} and \ref{def:reduced-problem-of-linear-elasticity}. If we assign
\begin{equation}
f=\left(F,\, \widetilde{T}\right),\qquad g={\bf C}\, {\rm E}(u_0),
\label{eq:loads_to_plug2}
\end{equation}
where  the operator ${\rm E}$ defined as \eqref{eq:E-def1}--\eqref{eq:E-def2} and $\widetilde{T}\in  [H^{-1/2}_{00}(\Gamma_{\rm T})]^n$ is given as functional
\[
\widetilde{T}: [H^{1/2}_{00}(\Gamma_{\rm T})]^n \to \mathbb{R},
\]
\begin{equation}
\widetilde{T}:u \mapsto \int\limits_{\Gamma_{\rm T}} T\cdot u\, dS,
\label{eq:widetilde-T-def}
\end{equation}
then equation \eqref{eq:linear_problem2_eq1} is equivalent to \eqref{eq:le-1}--\eqref{eq:le-3}, and equation \eqref{eq:linear_problem2_eq2} is equivalent to \eqref{eq:le-4}.
Moreover, linear mapping
\begin{gather*}
[L^2(\Gamma_{\rm T})]^n \to [H^{-1/2}_{00}(\Gamma_{\rm T})]^n,\\
T\mapsto \widetilde{T} 
\end{gather*}
is continuous.
\label{prop:linear_problem2_with_loads}
\end{Proposition}
\noindent {\bf Proof.} 
Equivalence between  \eqref{eq:linear_problem2_eq1} and \eqref{eq:le-1}--\eqref{eq:le-3} follows the lines of Proposition \ref{prop:linear_problem1_with_loads}.

We verify that $\widetilde{T}\in [H^{-1/2}_{00}(\Gamma_{\rm T})]^n$ when $T\in [L^2(\Gamma_{\rm T})]^n$. Since a generic element of $[H^{1/2}_{00}(\Gamma_{\rm T})]^n$ is $u|_{\Gamma_{\rm T}}$ for an arbitrary $u\in H^{1/2}_{\rm D}(\Gamma)$ (see Definition \ref{def:L-M-spaces}), we have 
\begin{align*}
\left|\left\la \widetilde{T}, u|_{\Gamma_{\rm T}} \right\ra\right| = \left|\,\int\limits_{\Gamma_{\rm T}} T\cdot u|_{\Gamma_{\rm T}}\, dx\right| &\leqslant \|T\|_{L^2}\, \|u|_{\Gamma_{\rm T}}\|_{[L^2(\Gamma_{\rm T})]^n}=\\
&=\|T\|_{L^2}\, \|u\|_{[L^2(\Gamma)]^n} \leqslant \|T\|_{L^2}\,\|u\|_{[H^{1/2}(\Gamma)]^n} = \|T\|_{L^2}\,\|u\|_{[H^{1/2}_{00}(\Gamma_{\rm T})]^n},
\end{align*}
where the last inequality follows from \eqref{eq:H1/2-norm-def} and the last equality follows from \eqref{eq:L-M-space-norm-def}. Using the same reasoning, we can verify that for $T_1, T_2\in[L^2(\Gamma_{\rm T})]^n$
\begin{align*}
\left\|\widetilde{T}_1-\widetilde{T}_2 \right\|_{[H^{-1/2}_{00}(\Gamma_{\rm T})]^n}&= \sup_{v\in[H^{1/2}_{00}(\Gamma_{\rm T})]^n\setminus\{0\}}\frac{\left|\left\la \widetilde{T}_1-\widetilde{T}_2, v \right\ra\right|}{\|v\|_{[H^{1/2}(\Gamma)]^n}}\leqslant
\\[1mm]& \leqslant \sup_{v\in[H^{1/2}_{00}(\Gamma_{\rm T})]^n\setminus\{0\}}\frac{\|T_1-T_2\|_{L^2}\,\|v\|_{L^2}}{\|v\|_{[H^{1/2}_{00}(\Gamma_{\rm T})]^n}}\leqslant \|T_1-T_2\|_{L^2},
\end{align*}
i.e. relation $T\mapsto \widetilde{T}$ is continuous.

Let $\sigma$ satisfy \eqref{eq:linear_problem2_eq2}. With $f=\left(F, \widetilde{T}\right)$ and  \eqref{eq:D1-def2},\,\eqref{eq:D1C-def2} plugged in, it becomes
\begin{equation}
\begin{pmatrix}-{\rm div}\, \sigma\\ {\rm tr}_{\nu, \Gamma_{\rm T}}(\sigma)\end{pmatrix}=\begin{pmatrix}F\\\widetilde{T}\end{pmatrix}.
\label{eq:div-and-trace-to-plug}
\end{equation}
By Green's formula \eqref{eq:Green's_formula_for_matrices2} we have that for any $u\in [W_{\rm D}]^n$
\[
\left \la {\rm tr}_{\nu, \Gamma_{\rm T}}(\sigma), {\rm tr\,}_{\Gamma_{\rm T}}(u) \right\ra = \intOm \tau: {\rm E} (u)\, dx+  \intOm ({\rm div}\, \sigma) \cdot u \, dx.
\]
By rearranging the terms, we get
\[
\intOm \sigma: {\rm E} (u)\, dx = \intOm (-{\rm div}\, \sigma) \cdot u \, dx + \left \la {\rm tr}_{\nu, \Gamma_{\rm T}}(\sigma), {\rm tr\,}_{\Gamma_{\rm T}}(u) \right\ra,
\]
and by plugging \eqref{eq:div-and-trace-to-plug} and $\widetilde{T}$, we find that
\[
\intOm \sigma: {\rm E} (u)\, dx = \intOm F \cdot u \, dx + \int\limits_{\Gamma_{\rm T}} T\cdot {\rm tr\,}_{\Gamma_{\rm T}}(u)\, dS,
\]
which is exactly \eqref{eq:le-4}.

Vice versa, let $\sigma\in L^{2}_{\rm sym}(\Omega)\cong L^2_{{\rm sym}, {\bf C}^{-1}}(\Omega)$ satisfy \eqref{eq:le-4}. In particular, one can choose $\delta u\in [W^{1,2}_{0}(\Omega)]^n\subset [W_{\rm D}]^n$ and then the last term in \eqref{eq:le-4} is zero. Hence, by Definition \ref{def:H_div_sym-def}, $\sigma\in H_{\rm sym}({\rm div}, \Omega)$  with $-{\rm div}\, \sigma = F$. From both \eqref{eq:le-4} and Green's formula \eqref{eq:Green's_formula_for_matrices2} we have that for any $u\in [W_{\rm D}]^n$
\[
\intOm F \cdot u \, dx + \int\limits_{\Gamma_{\rm T}} T\cdot {\rm tr\,}_{\Gamma_{\rm T}}(u)\, dS = \intOm \sigma: {\rm E} (u)\, dx  = \left \la {\rm tr}_{\nu, \Gamma_{\rm T}}(\sigma), {\rm tr\,}_{\Gamma_{\rm T}}(u) \right\ra +  \intOm (-{\rm div}\, \sigma) \cdot u \, dx.
\]
Since we already know that $-{\rm div}\, \sigma = F$, we can deduce that
\[
\int\limits_{\Gamma_{\rm T}} T\cdot {\rm tr\,}_{\Gamma_{\rm T}}(u)\, dS = \left \la {\rm tr}_{\nu, \Gamma_{\rm T}}(\sigma), {\rm tr\,}_{\Gamma_{\rm T}}(u) \right\ra.
\]
By the surjectivity of ${\rm tr}_{\Gamma_{\rm T}}$ (see Remark \ref{remark:restriction_of_traces_surjectivity}) we have \eqref{eq:div-and-trace-to-plug}, which, in turn, means \eqref{eq:linear_problem2_eq2}.
$\blacksquare$
\begin{Remark} Again, we note the continuity of linear mapping \eqref{eq:linear_solution2}--\eqref{eq:widetilde-T-def} acting as
\begin{gather*}
[W^{1,2}(\Omega)]^n \times [L^{2}(\Omega)]^n\times[L^{2}(\Gamma_{\rm T})]^n \to L^{2}_{\rm sym}(\Omega)\cong L^2_{{\rm sym}, {\bf C}^{-1}}(\Omega),
\\
(u_0, F, T)\mapsto \sigma.
\end{gather*}
Moreover, Theorem \ref{th:linear_problem2_solution} and Proposition \ref{prop:linear_problem2_with_loads} mean that the schematic representation of Fig. \ref{fig:spaces} is suitable for the static stress problem of  linear elasticity with the mixed boundary conditions.
\end{Remark}

\begin{Remark} Similarly to Remark \ref{remark:explicit-representations-displacement-bc}, for the problem with the mixed boundary conditions, space $\mathcal{V}$ and operators $P_{\,\mathcal{V}}$ and $Q$ can be expressed as the following.

\noindent \underline{Space $\mathcal{V}$ written explicitly.} From \eqref{eq:D1-def1}--\eqref{eq:D1-def2} we get that
\[
\mathcal{V} = {\rm Ker}\, {\bf D_{1, C}}  =\left\{\sigma\in H_{\rm sym}({\rm div}, \Omega): {\rm div}\, \sigma = 0 \text{ and } {\rm tr}_{\nu, \Gamma_{\rm T}}(\sigma)=0\right\}.
\]

\noindent \underline{Operator $P_{\,\mathcal{V}}$ written explicitly.} Let 
\[
P_{\,\mathcal{V}}\, g=\sigma
\]
for some $g, \sigma\in L^{2}_{{\rm sym}, {\bf C}^{-1}}(\Omega)$, which means  \eqref{eq:projection_abstract_def} with inner product \eqref{eq:weighted-ip}, i.e.
\[
\begin{cases} {\bf D_{1, C}},\sigma=0, \\[1mm]
\intOm\tau: {\bf C}^{-1}\sigma\, dx =\intOm \tau:{\bf C}^{-1} g\, dx & \text{for any }\tau \in L^{2}_{{\rm sym}, {\bf C}^{-1}}(\Omega)\text{ such that } {\bf D_{1, C}}\, \tau=0,
\end{cases}
\]
or
\[
\begin{cases} 
\sigma\in H_{\rm sym}({\rm div}, \Omega),\\[1mm]
{\rm div}\,\sigma=0, \\[1mm]
{\rm tr}_{\nu, \Gamma_{\rm T}}(\sigma)=0,\\[1mm]
\intOm\tau: {\bf C}^{-1}\sigma\, dx =\intOm \tau:{\bf C}^{-1} g\, dx & \text{for any }\tau \in H_{\rm sym}({\rm div}, \Omega) \text{ with } {\rm div}\, \tau=0\, \text{ and } \, {\rm tr}_{\nu, \Gamma_{\rm T}}(\tau)=0.
\end{cases}
\]

\noindent \underline{Operator $Q$ written explicitly.} Let
\[
Q\,f=\sigma
\]
for some $f\in [L^2(\Omega)]^n \times [H^{-1/2}_{00}(\Gamma_{\rm T})]^n, \,\sigma\in  L^{2}_{{\rm sym}, {\bf C}^{-1}}(\Omega)$. By the construction of $Q$ it means that
\[
\begin{cases}
\sigma\in {\rm Im}\, {\bf E_{1,C}} = \left({\rm Ker}\, {\bf D_{1, C}}\right)^{\perp_{{\bf C}^{-1}}},\\
{\bf D_{1, C}}\,\sigma=f,
\end{cases}
\]
i.e.
\[
\begin{cases}
\intOm \tau: {\bf C}^{-1}\sigma=0& \text{for any }\tau \in H_{\rm sym}({\rm div}, \Omega) \text{ with } {\rm div}\, \tau=0\, \text{ and } \, {\rm tr}_{\nu, \Gamma_{\rm T}}(\tau)=0\\[1mm]
\sigma\in H_{\rm sym}({\rm div}, \Omega),\\[1mm]
\begin{pmatrix}
-{\rm div}\, \sigma \\
{\rm tr}_{\nu, \Gamma_{\rm T}}(\sigma)
\end{pmatrix}
=f.
\end{cases}
\]
\end{Remark}

\section{Conclusions}
\label{sect:conclusions}
We have analyzed the static stress problems of linear elasticity with the mixed boundary conditions and with the displacement boundary conditions by representing them through the pairs of adjoint unbounded linear operators. Such tools of linear functional analysis prove to be sufficient  to show existence and uniqueness of the solution and its continuous and linear dependence on given loads.  In the future, the formulations of the current paper can be included in various nonlinear models by including phenomena of plasticity, friction or electromagnetism. It can also help to  construct conforming finite element approximations and prove guaranteed a posteriori error estimates and the two-sided energy estimates, see e.g. \cite{Krizek1996}, \cite{Necas1981}.

While we acknowledge, that our work is parallel to the classical variational theory in linear elasticity, we believe that our approach and presentation fills in an existing gap. In particular, a rigorous identification of the adjoint operator to the explicitly written stress-to-force mapping has not been done before for the case of mixed boundary conditions. Due to our focus on the theory of adjoint operators, we hope that the reader's intuition about the Fundamental Theorem of Linear Algebra could be carried over to problems in continuum mechanics. 

As we mentioned above, quite complex models in mechanics often include linear elasticity as a building block. For example, stress distributions in a model of elastoplasticty with hardening (see e.g. \cite[Section 13.5]{Necas1981}, \cite{Han2012}) also take values in the space $L^2_{{\rm sym}, {\bf C}^{-1}}(\Omega)$, but that is an evolution problem, which involves the intersection of the subspace $\mathcal{V}$
 with a closed convex (but not affine) set $\mathcal{P}-\sigma(t)$, where  $\sigma(t)$ is the stress variable we found in the current paper, cf. Fig. \ref{fig:spaces}. Such models will be the topic of further research.


\appendix
\section{A connection to the variational formulations and the displacement problems}
\label{sect:appendix-on-other-problems}
\subsection{Variational formulation for the stress problems}
We would like to clarify the connection between the formulations we derived earlier and the classical variational formulations. Specifically, in the abstract setting of Proposition \ref{prop:inverse_on_U} and Theorem \ref{th:abstract_linear_framework_solution} we can consider the following minimization problem in a Hilbert space $\mathcal{H}$:
\begin{equation}
\underset{\sigma\in \mathcal{V}+Qf}{\rm arg\, min}\, \mathscr{S}(\sigma), \qquad
\mathscr{S}(\sigma) = \frac{1}{2} \la\sigma, \sigma\ra - \la g, \sigma\ra,
\label{eq:abstract_variational}
\end{equation}
It can be shown that the minimum point coincides with the solution $\sigma$ in Theorem \ref{th:abstract_linear_framework_solution}. 

By plugging in \eqref{eq:abstract_variational}  the spaces and parameters from Section \ref{ssect:linear_problem1_solution} or Section \ref{ssect:linear_problem2_solution} we obtain the variational formulations of the respective problems in terms of complementary energy $\mathscr{S}$. In the latter case (the problem with the mixed boundary conditions) we obtain the formulation of \cite[(5.6)-(5.7), p.~106]{Necas1981} exactly. Notice, that with such substitutions the Hilbert space $\mathcal{H}$ is $ L^{2}_{{\rm sym},{\bf C}^{-1}}(\Omega)$ and the inner product in \eqref{eq:abstract_variational} is given by \eqref{eq:weighted-norm}.

\subsection{Operator formulation for the displacement problems}
To recover the full solution $(u, \varepsilon, \sigma)$ in the sense of Definitions \ref{def:static-problem-of-linear-elasticity-mixed}--\ref{def:static-problem-of-linear-elasticity-displacement} from the stress solution $\sigma$ one can set
\[
\varepsilon = {\bf C}^{-1} \sigma,
\]
and then use Proposition \ref{prop:inverse_on_U} with
\[
\mathcal{H} = [L^{2}(\Omega)]^n,\qquad \mathcal{U} = [L^{2}(\Omega)]^n, \qquad \mathcal{V} = \{0\}, \qquad \mathcal{X}= L^{2}_{\rm sym}(\Omega), \qquad A= {\bf E_{0}}
\]
(see Section \ref{ssect:aux_operators})  for the static problem of linear elasticity with the displacement boundary conditions, or with 
\[
\mathcal{H} = [L^2(\Omega)]^n\times [H^{1/2}_{00}(\Gamma_{\rm T})]^n, \qquad \mathcal{U} = [L^2(\Omega)]^n\times [H^{1/2}_{00}(\Gamma_{\rm T})]^n, \qquad \mathcal{V} = \{0\},
\]
\[
\mathcal{X}= L^{2}_{\rm sym}(\Omega), \qquad A= {\bf E_{1}}
\]
(see Section \ref{ssect:main_operators}) for the static problem of linear elasticity with the mixed boundary conditions. Recall that both operators ${\bf E_{0}}$ and ${\bf E_{1}}$ are injective (see Corollaries \ref{cor:E0-D0-properties} and \ref{cor:E1-D1-properties}, respectively), therefore we can, indeed, apply Proposition \ref{prop:inverse_on_U} with $\mathcal{V} = \{0\}$.
\subsection{Variational formulation for the displacement problems}
Variational formulations in terms of the displacements variables are well-known, see e.g. \cite[Remark 2.7, pp.~94--95]{Necas1981}. For completeness, we will provide them using our notation. For the problem of Definition \ref{def:static-problem-of-linear-elasticity-displacement} (with the displacement boundary conditions) we write
\begin{equation}
\underset{u\in [W_0^{1,2}(\Omega)]^n+u_0}{\rm arg\, min}\, \Phi(u), \qquad
\Phi(u) = \frac{1}{2} \intOm {\rm E}(u): {\bf C}{\rm E}(u)\, dx - \intOm F\cdot u\, dx.
\label{eq:variational_disp1}
\end{equation}
For the problem of Definition \ref{def:static-problem-of-linear-elasticity-mixed} (with the mixed boundary contitions) we write
\begin{equation}
\underset{u\in [W_D]^n+u_0}{\rm arg\, min}\, \Phi(u), \qquad
\Phi(u) = \frac{1}{2} \intOm {\rm E}(u): {\bf C}{\rm E}(u)\, dx - \intOm F\cdot u\, dx - \int\limits_{\Gamma_{\rm T}} T\cdot {\rm tr}_{\Gamma_T}(u)\, dS.
\label{eq:variational_disp2}
\end{equation}
Unlike the case of the stress problems, our framework of unbounded operators of Sections \ref{ssect:aux_operators}, \ref{ssect:main_operators} is inconvenient for such formulations, while operator ${\rm E}$ is bounded (Proposition \ref{prop:E-is-bounded}), so that the problems \eqref{eq:variational_disp1}, \eqref{eq:variational_disp2} admit a solution via the classical Lax-Milgram lemma (see e.g.  \cite[Corollary 5.8, p.~140]{Brezis2011}, \cite[Section 1.1, pp.~1--10]{Gatica2014}, \cite[Section 5.1, pp.~89--93]{Sayas2019}).

\section*{Acknowledgements}
The authors thank Pavel Krej\v{c}\'{i} and Giselle Antunes Monteiro from the Institute of Mathematics CAS for helpful scientific discussions.
\newline\noindent{\bf Funding:} both authors were supported by the GA\v{C}R project GA24-10586S and the Czech Academy of Sciences (RVO: 67985840).
\section*{Statements and Declarations}
\noindent{\bf Conflict of Interest:} the authors declare that they have no conflict of interest. \newline
\noindent{\bf Data availability:} no datasets were generated or analyzed during the current study.

\bibliographystyle{plainnat} 
\bibliography{linear_elasticity-R1.bib}





\end{document}